\documentclass[a4paper,oneside,11pt]{article}
\usepackage{amsfonts}
\usepackage{amsmath,amssymb}
\usepackage{amssymb}
\usepackage{makeidx}
\usepackage{graphicx}
\usepackage{float}
\usepackage{amsthm}
\usepackage[hmargin=3cm,vmargin=3.5cm]{geometry}
\theoremstyle{plain}
\newtheorem{counter}{Counter}[subsection]
\newtheorem{lema}[counter]{\scshape{Lemma}}
\newtheorem{propozitie}[counter]{\scshape{Proposition}}
\theoremstyle{definition}
\newtheorem{definitie}[counter]{\scshape{Definition}}
\newtheorem{conventie}[counter]{\scshape{Convention}}
\theoremstyle{remark}
\newtheorem*{remarca}{Remark}
\newtheorem*{exemplu}{Example}

\begin{document}
\author{Dan-Titus Salajan}
\pagenumbering{roman}
\title{$CAT(0)$ geometry for The Thompson Group}
\date{\today}
\maketitle
\cleardoublepage
\newpage
\begin{center}
\bf Abstract
\end{center}
We investigate Farley's $CAT(0)$ cubical model for Thompson's group $F$ (we adopt the classical language of $F$, using binary trees and piecewise linear maps). Main results include: in general, Thompson's group elements are parabolic; we find simple, exact formulas for the $CAT(0)$ translation lengths, in particular the elements of $F$ are ballistic and uniformly bounded away from zero; there exist flats of any dimension and we construct explicitly many of them; we reveal large regions in the Tits Boundary, for example the positive part of a non-separable Hilbert sphere , but also more complicated objects. En route, we solve several open problems proposed in Farley's papers.

\vspace{3mm}


\clearpage

\clearpage

\begin{center}
$\clubsuit$
\end{center}

\clearpage

\tableofcontents

\clearpage

\begin{center}
$\clubsuit$
\end{center}

\clearpage

\pagenumbering{arabic}

\section{Introduction}

In this chapter we set up the scene of this work and discuss the main results.

\subsection{History and Results}

Thompson's group $F$ is the group of all piecewise linear homeomorphisms of the unit interval, with finitely many slopes, all breaking points located at dyadic rational numbers and all the derivatives are powers of $2$. It was discovered by Richard Thompson in the 60's and rediscovered by topologists in late 70's. For a light introduction to $F$ we recommend \cite{notices} and for a technical introduction \cite{intro}. $F$ is finitely presented, it doesn't contain non-abelian free groups \cite{free}, while its amenability status is unknown and legendary (to sense a tension we cannot describe here, we refer to  \cite{funny1}). For a long list of unsolved problems about the Thompson group one can look at \cite{probleme1,probleme2}.

In this thesis we study Thompson's group from a $CAT(0)$ point of view. In a sequence of papers \cite{farley03, farley05,farley08}, Farley introduced and studied $CAT(0)$ cubical models for the class of diagram groups \cite{diagrams} which includes the Thompson group (see the next two sections for details about the model). As it is acknowledged in \cite{farley03}, in the case of $F$, the model was studied before (without the $CAT(0)$ condition) by Stein \cite{stein}, and even before in a disguised form \cite{higman}. The formal language we adopt is the one in \cite{stein} with the only difference that the multiplications are reversed for convenient specification (see Section 1.3 for the precise definition of the complex). The complex is locally finite, but infinite dimensional, the action of $F$  is proper but not cocompact. As a consequence, Farley proves the Haagerup property for $F$ (in fact, the induced action on the hyperplanes of the complex gives an equivariant Hilbert embedding with compression $\sqrt{n}$, and any asymptotical improvement (no matter how small) of this bound would imply amenability of $F$; however there is not much room to improve- indeed the upper bound is $\sqrt{n}log(n)$ \cite{ags}).

While in our work we focus on $CAT(0)$ geometry questions (related with $F$), it is worthwhile to point out some interesting facts about possible applications of the $CAT(0)$ theory to group theoretical questions (about $F$). Indeed, since $F$ is not a $CAT(0)$ cubical group, but still acts naturally on a $CAT(0)$ cubical complex, it is interesting to see what properties do they share.(The large class of) $CAT(0)$ cubical groups have quadratic Dehn function, are automatic, can be embedded with high Hilbert compression and have Yu's Property $A$ \cite{niblo1,niblo2,book}. All these properties are open for the Thompson group with the notable exception of the Dehn function which was computed in \cite{gubadehn} (which is also a first step toward deciding if $F$ is automatic). It looks challenging to extend these techniques to the infinite dimensional setting. We should also note that Farley's motivation to study this complex was deciding amenability of $F$ via Adams-Ballman Theorem. This strategy fails since he finds a Tits arc of length $\pi/2$ consisting of globally fixed points, and indeed, it is shown in \cite{cm} that this strategy could never work. However, Sageev asks on the Geometric Group Theory Problems Wiki for more subtle $CAT(0)$ strategies to address this question.

Our first theorem points out a modest difference with the classical finite-dimensional setting. By a theorem of Bridson, the isometries of a finite-dimensional $CAT(0)$ cubical complex are semi-simple. Taking this fact as evidence and also group theoretical results of Guba and Sapir, Farley conjectured that all the isometries of the complex are semi-simple, Question 3.17 in \cite{farley03}. We show that the situation is more complicated and solve this problem in the negative. We call an element of $F$ \emph{irreducible} if it fixes only $0$ and $1$ (as a homeomorphism of the interval). We denote by $X$ the associated $CAT(0)$ cubical complex.

\vspace{3mm}

\noindent \textbf{Theorem 1.} \emph{ Let $g \in F$ be an irreducible element. Then the $CAT(0)$ translation length of $g$ in $X$ equals $\sqrt{log_2(g'(0))^2+log_2(g'(1))^2}$ and, if $g'(0)g'(1) \neq 1$ then $g$ is parabolic.}

\vspace{3mm}

Since natural examples of groups acting on $CAT(0)$ cubical complexes by parabolic isometries are quite rare, we consider this situation quite interesting. Also, a feature of this theorem might be the exact computation of the translation lengths, which is done by using medium scale geometry arguments. A consequence of the proof is also that all the non-identity elements of $F$ are ballistic (positive translation lengths) and the translation lengths are uniformly bounded away from zero.

\vspace{3mm}

\noindent \textbf{Corollary 1.}\emph{ All non-identity elements of $F$ are ballistic with translation lengths uniformly bounded away from zero. More precisely any translation length is at least $\sqrt{2}$, the constant being sharp.}

\vspace{3mm}

Some elements of $F$ are hyperbolic in the model (it will appear clear in the text that they are combinatorial coincidences). Using them we can explicitly construct and locate flats of any dimension.

\vspace{3mm}

\noindent \textbf{ Theorem 2.}\emph{ $X$ contains flats of any dimension.}

\vspace{3mm}

We now discuss facts connected with the Tits Boundary. While the flats and the ballistic isometries induce rich geometry at the infinity of $X$, we clarify an independent part of the boundary corresponding to a remarkable sub-complex which is a model for the geometry of finite, rooted, ordered, binary trees (or equivalently for the dyadic partitions of the unit interval).

To describe this large region of the boundary we need a definition. Denote by $I_{n,k}$ the standard dyadic intervals $[\frac{k-1}{2^n}, \frac{k}{2^n}]$, where $n \geq 1$ is a positive integer and $k=1, \ldots, 2^n$. These intervals are in natural one-to-one correspondence with the vertices of the full rooted, ordered binary tree. A map from the set of all standard dyadic intervals (or equivalently from the full binary tree) to $[0,1]$ is called \emph{a flow} if it satisfies the following two properties: $V([0,1])=1$ and $V(I_{n,k})=V(I_{n+1,2k-1})+V(I_{n+1,2k})$ for any $n$ and $k$.

\vspace{3mm}

\noindent \textbf{Theorem 3.} \emph{The Tits boundary of $X$ contains a copy of the metric space $(M, \rho)$, where $M$ is the set of all flows and $\rho$ (giving the Tits angle) is defined, for any two flows $V$ and $W$, by the formula:}

$$\rho(V,W)=arccos {(lim_{n \rightarrow \infty}\sum_{k=1}^{2^n}\sqrt{V(I_{n,k})W(I_{n,k})})}$$

\emph{This region is invariant under $F$ and the action is defined in the following way: notice first that specifying the values of a flow at all but finitely many intervals still determines uniquely the flow; notice also that any element of $F$ maps linearly all the standard dyadic intervals to standard dyadic intervals, with finitely many exceptions. The action is given by $(gV)(P)=V(g(P))$ for any $g \in F$, any flow $V$ and any standard dyadic interval $P$ which is mapped linearly to another standard dyadic interval by $g$.}

\vspace{3mm}

To have a feeling of how huge is this (small) region of the boundary, we highlight a minuscule part of it (for details see Section 3.1).

\vspace{3mm}

\noindent \textbf{Corollary 2.} \emph{The Tits boundary of $X$ contains the positive orthant of a non-separable Hilbert sphere with the angular metric.}

\vspace{3mm}

Much remains to be understood about the Tits boundary, for example the interaction of the spheres coming from the flats with this region, but also the behavior of the canonical points associated with the elements of $F$. However, we can push one more remarkable conclusion. Farley \cite{farley08} proposed an interesting combinatorial approximation of the Tits Boundary of a cube complex. A profile is roughly a collection of hyperplanes which are likely to be crossed by a geodesic ray. He divides the boundary on such classes of profiles. In his search for global fixed points Farley analyzes all the profiles, except (a very large) one. A globally fixed arc of length $\pi/2$ is found and any other fixed point should lie in that remaining profile. Combined with Farley's work, our Theorem 3 is able to eliminate this possibility and solves the main question left open in \cite{farley08} (Conjecture 7.8).

\vspace{3mm}

\noindent \textbf{Corollary 3.} \emph{ The Thompson group fixes at infinity of $X$ a Tits arc of length $\pi/2$ and no other point.}

\vspace{3mm}

We end up by saying that while the profiles offer a very helpful guide at infinity for our space, it is not always a reliable method. Indeed, Farley conjectured (Conjecture 2.8(1)) in \cite{farley08} that for any profile of a locally finite $CAT(0)$ cubical complex there is a geodesic ray realizing it. In the Appendix, we construct a simple complex with one point at infinity and two profiles, contradicting Farley's proposal.

\vspace{3mm}

\noindent \textbf{Overview.}

\vspace{2mm}

In the rest of this first chapter, we describe the model and its basic properties. Section 1.4. contains important remarks. In the second chapter we prove Theorem 1 and Theorem 2. In the third chapter we prove Theorem 3. More detailed descriptions are given at the beginning of each chapter.

A few words about the necessary background needed to read this document. About Thompson's group we really only use its description with pairs of binary trees. However, the setting of this thesis is $CAT(0)$ geometry and cubical complexes. We use freely the most basic $CAT(0)$ theory \cite{book}: definition, angles, projections, boundary. About cubical complexes we also use freely the most basic theory \cite{sageev,chepoi} : hyperplanes and the conditions on the link (which have been now established in maximal generality in the appendix of \cite{leary}). All the rest is defined. Also, we believe to have a good system of specification of binary trees and associated diagrams. Since most of the time the graphical Thompson-like computations we do are trivial, we just claim well defined equalities between diagrams (unfortunately, sometimes one has to draw a small pictures to check it, but this is the specificity of the subject).

\subsection{Notation}

We formally introduce the necessary objects to define our model and the unpleasant associated notations (but which best describes their simple graphical visualization).

\begin{definitie}
Let $n \geq 1$ be a positive integer. A function $g: [0,1] \rightarrow [0,n]$ is called \emph{a Thompson-like function} (of degree $n$) if it is a piecewise linear homeomorphism with finitely many breaking points, all located at dyadic rationals, and with all the derivatives powers of $2$. Notice
that such a function maps dyadic rationals to dyadic rationals. If $n=1$ then $g$ is an element of the Thompson Group $F$. For any $n \geq 1$ we denote by $F_n$ the set of all Thompson-like functions of degree $n$ and by $F_{\infty}=\bigcup_{n=1}^{\infty}F_n$ the set of all Thompson-like functions.
\end{definitie}

Like in the case of the Thompson group \cite{intro}, we have an alternative description of Thompson-like functions using tree diagrams. We first introduce some general notations for binary trees and then state some lemmas which should be clear for the $F$-familiar reader.

\begin{conventie}
A finite, rooted, ordered, binary tree will be simply called \emph{tree}. Also, if $s \geq 1$ is an integer, $[s]$ will always denote the set $\{1, \ldots, s\}$.
\end{conventie}

Before giving a formal definition, let us notice that any tree is made out of a finite number of \emph{carets} (a caret is a vertex, called \emph{the head of the caret}, together with its two children and with the corresponding edges). The number of carets of a tree is called \emph{the size} of the tree and we denote this quantity by $|T|$. It is clear (and easy to prove by induction) that a tree of size $n$ has $2n+1$ vertices, out of which $n+1$ are \emph{leaves}, that is vertices without children. The \emph{empty tree} is the unique tree of size zero and will be denoted $*$. The unique tree of size one will be denoted $\wedge$.

For a given tree $T$, we call (a) \emph{a free caret} of $T$, a caret with both its children leaves (b) \emph{blocked caret} a caret with none of its children leaves and (c) \emph{mixed caret} a caret with one children leaf and the other one not a leaf. It is funny to notice that for any non-empty tree we have that the number of free carets is one plus the number of blocked carets (easy to prove by induction on the number of carets). We denote by $||T||$ the number of free carets of $T$.

Let us give precise definitions now. It will be convenient to specify a tree as follows: we first index the positions of the carets in the full (infinite) rooted binary tree and then a tree will be just a set of positions (respecting a gluing rule). In the following definition, $(k,l)$ will interpret the $l^{th}$ caret (from left to right) at depth $k$.

\begin{definitie}
\label{def tree}

 For a subset $I$ of the natural numbers we denote by $2I$ the set $\{2i-1,2i | i \in I\}$. The positions of the carets are specified by the sets $Bin^k=\{k\} \times [2^k]$ (for any natural $k$) and $Bin=\bigsqcup_{k=0}^{\infty}Bin^k$. A tree $T$ is a finite subset of $Bin$ with the following gluing condition: $T^k \subseteq 2T^{k-1}$ for all $k$, where $T^j=p_2(Bin^j \cap T)$ and $p_2$ is the second projection of $Bin$. An infinite tree is defined the same but, of course, without the finiteness condition. For example the empty tree will be the empty set, the full binary tree will be $Bin$ and if we denote by $T_n$ the full tree of depth $n$, then $T_3=\{(0,1),(1,1),(1,2), (2,1),(2,2),(2,3),(2,4)\}$. Sometimes we will refer to a caret in a given tree by specifying its position. We denote by $Trees$ the set of all trees. It is easy to see that for any two (infinite) trees their union and intersection are also ( possible infinite) trees.
\end{definitie}

\begin{conventie}
\label{conventie lipire}
For a tree of size $n$ we always index its leaves from left to right with $1,2, \ldots, n+1$. Another important notation is the following (gluing trees): given a tree $A$ of size $n$, $I \subseteq [n+1]$ and $B_1, \ldots, B_{|I|}$ some trees, we denote $A \oplus^I (B_1, \ldots, B_{|I|})$ the tree obtained by simultaneously attaching each $B_i$ at the leaf of $A$ indexed by the $i^{th}$ element of $I$. For simplicity, if all $B_i$'s are the (unique) tree of size one (i.e. just one caret) we just denote $A \oplus^I$. If $|I|=1$ (i.e. we glue just one tree) we just write $A \oplus^i B$ (in particular if the size of $B$ is one, we denote $A \oplus^i$). When we attach trees at all the leaves we simply write $A \oplus (T_1, \ldots, T_{n+1})$ instead of $A \oplus^{[n+1]} (T_1, \ldots, T_{n+1})$. In particular, when we attach a caret to all the leaves, we simply write $T \oplus$.  \end{conventie}

\begin{exemplu}
The previous notations allow very useful specification of trees. For example the most left tree of depth $n$ denoted $L_n$, can be defined recursively as follows: $L_0=*$ and $L_k=L_{k-1}\oplus^1$. The most right tree of
depth $n$, denoted $R_n$, can be defined similarly: $R_0=*$ and $R_k=R_{k-1} \oplus^k$. The full tree of depth $n$, denoted $T_n$, can be defined: $T_0=*$ and $T_k=T_{k-1}\oplus^{[2^{k-1}]}$. Also, with our notation, the two standard generators of the Thompson Group $F$ are $(L_2,R_2)$ and $(L_3,L_2 \oplus^2)$ (see \cite{intro}).
\end{exemplu}

\begin{definitie}
A \emph{dyadic rational} is a rational number of the form $k/2^n$ where $k$ and $n$ are integers. A \emph{standard interval} is an interval of the form $[k/2^n,(k+1)/2^n]$, where $k$ and $n$ are nonnegative integers. For $n \geq 1$ we call \emph{a standard partition} of $[0,n]$ a finite sequence $0=x_0 < x_1 < \ldots < x_s=n$ such
that each interval $[x_i,x_{i+1}]$ is a standard interval.
\end{definitie}

The proof of the next lemma is easy and can be found in \cite{intro}.

\begin{lema}
There is a bijection between the standard partitions of $[0,1]$ and the set of all trees, such that each leaf of the tree correspond to an interval partition in the obvious way: the interval $[\frac{k-1}{2^n}. \frac{k}{2^n}]$ correspond to the $k^{th}$ vertex at depth $n$ (in the full binary tree).
\end{lema}

Similarly, we have:
\begin{lema}
There is a bijection between the standard partitions of $[0,n]$ and the set of all $n$-tuples of trees.
\end{lema}
\noindent \emph{Proof.} Notice first, that a standard partition of $[0,n]$ contains all the integers $0,1, \ldots, n$.  Also, the set of all standard partitions of $[j,j+1]$ is in bijection with the set of all standard partitions of $[0,1]$ (via the map "substracting $j$"). With the previous Lemma we are done. $\clubsuit$

\vspace{2mm}

\begin{lema}
Let $g:[0,1] \rightarrow [0,n]$ be a Thompson-like function. Then there is a standard partition of $[0,1]$:

$$x_0=0 < x_1 < x_2 < \ldots < x_s=1$$

\noindent such that $g$ is linear on each interval $[x_i,x_{i+1}]$ and
$$0=f(x_0)< f(x_1) < \ldots < f(x_s)=n$$
is a standard  partition of $[0,n]$.
\end{lema}

\noindent \emph{Proof.} Warning: The symbol $\bigsqcup$ denotes, only in this proof, a union with disjoint interiors.

Let $0=x_1 < x_2< \ldots <x_r=1$ be the breaking points of $g$. If we manage to write every interval  $[x_i,x_{i+1}]$ as a disjoint
union of elementary dyadic intervals such that the images are also elementary intervals then we are done. So fix an interval of the partition
$[x_i,x_{i+1}]$. If $f$ has slope $2^p$, $p \in \mathbb{Z}$, then choose a very large $m>0$ (much larger than $p$) such that we can write
$x_i=k/2^{m+p}$, $x_{i+1}=l/2^{m+p}$, $f(x_i)=k'/2^m$ and $f(x_{i+1})=l'/2^m$. Notice that we have to have $l-k=l'-k'$ in this case. Then we can
write $[x_i,x_{i+1}]$ as being the partition $$[k/2^{m+p},(k+1)/2^{m+p}] \sqcup [(k+1)/2^{m+p},(k+2)/2^{m+p}] \sqcup \ldots \sqcup [(l-1)/2^{m+p},l/2^{m+p}]$$
\noindent which maps linearly to the partition $$[k'/2^{m},(k'+1)/2^{m}] \sqcup [(k'+1)/2^{m},(k'+2)/2^{m}] \sqcup \ldots \sqcup
[(l'-1)/2^{m},l'/2^{m}].$$ $\clubsuit$

\vspace{2mm}

\begin{definitie}
An $n$-diagram is an ordered list of $n+1$ trees $(T,T_1, \ldots, T_n)$ such that the number of leaves in $T$ equals the total number of leaves in the list $(T_1, \ldots, T_n)$. The leaves of $T_1, \ldots, T_n$ are indexed continuously (left to right) from $1$ to $|T|+1$.
\end{definitie}

Notice that given an $n$-diagram there is a unique $n$-degree Thompson-like function $g$ such that $g$ maps linearly the interval indexed by the $k^{th}$ left leaf of the diagram to the interval indexed by the $k^{th}$ right leaf. Conversely, by Lemma 1.2.8, to each $n$-degree Thompson-like function $g$ we can attach a $n$-diagram with the same property.

In fact, we can attach many diagrams to a given Thompson-like function $g$. Fix such a diagram representing $g$. We can construct new diagrams corresponding to $g$  by the following two operations:

(a) \emph{adding a caret}: by attaching a caret to the $k^{th}$ leaf,
both on the left and on the right. The new diagram still represents $g$.

(b) \emph{removing a caret}: if for some $k$, for both trees of the diagram the leaves $k$ and $k+1$ belong to a common caret, one can remove the two carets. The new diagram still represents $g$.

\begin{definitie}
Two diagrams are called \emph{equivalent} if one can be obtained from the other by a sequence of operations (adding or removing a caret). This is clearly an equivalence relation. A diagram is called \emph{reduced} if one cannot remove a common caret. In each equivalence class of diagrams, there is exactly one reduced diagram.
\end{definitie}

\begin{lema}
There is a bijection between Thompson like-functions and classes of diagrams.
\end{lema}

\vspace{2mm}

\noindent \emph{Proof.} We claim that the function that maps a diagram class to its corresponding Thompson-like function is a bijection. It is a surjection by the previous discussion. We only have to prove that two distinct reduced diagrams correspond to two different functions. Add enough common carets to get diagrams $A'$ and $B'$ equivalent with $A$ and $B$ and with the same tree in the first position. Let $k_0$ be the smallest positive integer
such that $A'$ and $B'$ differ at the leaf $k_0$. Let $f$ and $g$ be the functions corresponding to the diagrams $A'$ and $B'$. Then it follows that the slope of $f$ and $g$ are different at $k_0$. $\clubsuit$

\begin{figure}[!h]
\label{rolex}
\begin{center}
\includegraphics{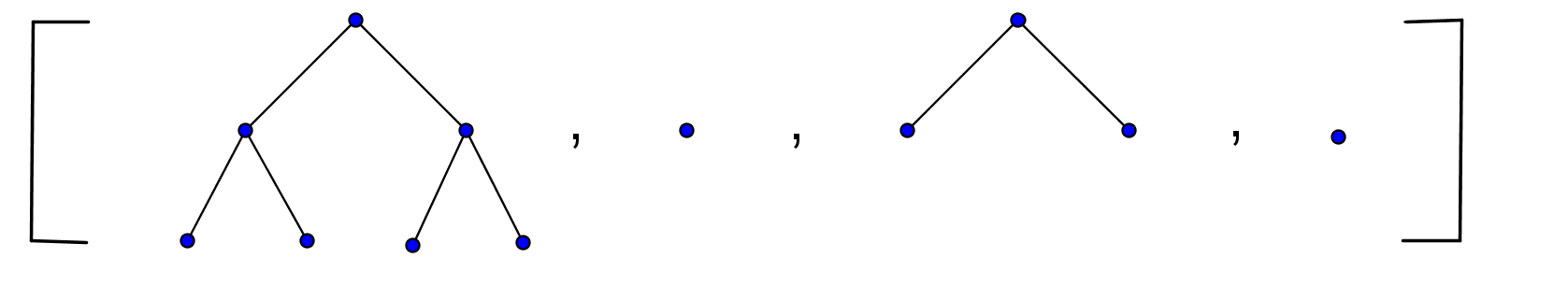}
\caption{A Thompson-like function $[0,1] \rightarrow [0,3]$ which maps linearly $[0,1/4]$ onto $[0,1]$, $[1/4,2/4]$ onto $[1,3/2]$, $[2/4,3/4]$ onto $[3/2,2]$ and $[3/4,1]$ onto $[2,3]$.}
\end{center}
\end{figure}

\begin{conventie}
Like in the case of $F$ there is a simple way to compose functions using diagrams. If $f \in F$ and $g$ is a Thompson-like function of degree $n$, then taking any two diagrams representing $f$ and $g$, with sufficiently many carets added, say of the form $(T,S)$ and $(S,T_1, \ldots,T_n)$, then it is easy to see that $(T,T_1, \ldots, T_n)$ is a diagram representing $g \circ f$. This is independent of the choice of the diagrams (of this form). The obvious right action of $F$ on the set of all Thompson-like function becomes a left one if we always consider the multiplication in the order of the graphical calculus (both in $F$ and when we work with Thompson-like functions). Since diagrams are more ubiquitous in this thesis, we prefer using them and acting on the left.
\end{conventie}

\subsection{The Model}

In this section we describe the model. We need a couple of additional definitions. Our language is very close to the one in \cite{stein}.

\begin{definitie}
Let $T$ be a non-empty tree. The \emph{left side} and \emph{right side} of $T$ are the unique trees, denoted $T^l$ and $T^r$ with the property that $T= \wedge\oplus^{\{1,2\}}(T^l,T^r)$, i.e. the two trees obtained by removing the top caret of $T$.
\end{definitie}

\begin{definitie}
Let $\Delta=(T,T_1, \ldots, T_n)$ be a diagram and $j \in [n]$. We denote by $\Delta(j)$ the diagram $$(T,T_1, \ldots,T_{j-1},T_j^l,T_j^r,T_{j+1}, \ldots, T_n)$$ if $T_j \neq *$ and $$(T \oplus^k,T_1, \ldots,T_{j-1},*,*,T_{j+1}, \ldots, T_n)$$ if $T_j=*$, where $k$ is the index of the unique leaf of $T_j$.

The process of obtaining $\Delta(j)$ from $\Delta$ is called \emph{cutting the diagram $\Delta$ at $j$} and the reverse process of getting $\Delta$ from $\Delta(j)$ is called \emph{gluing at the positions $(j,j+1)$}.
\end{definitie}

\begin{remarca}
Notice that if two diagrams $\Delta$ and $\Delta'$ represent the same Thompson-like function $g$, then $\Delta(j)$ and $\Delta'(j)$ will also represent the same Thompson-like function, namely the function obtained from $g$ by doubling all the slopes in the interval $g^{-1}([j-1,j])$ leaving the function unchanged on the interval $[0,j-1]$ and shifting the image from $[j,n]$ to $[j+1,n+1]$. If $g$ is a (degree $n$) Thompson-like function we denote by $g[j]$ the well defined (degree $n+1$) Thompson-like function obtained by cutting at $j$. More formally, $g[j]=s_{n,j} \circ g$, where $s_{n,j}$ is the unique piecewise linear bijection from $[0,n]$ to $[0,n+1]$, with all the slopes $1$, except on the interval $[j-1,j]$ where the slope is $2$.
\end{remarca}

The vertices of the cubical complex will be the set of all Thompson-like functions and the underlying graph structure will be given by the cutting/gluing operation: $f,g$ (or their associated diagrams) will be at distance $1$ if and only if $g$ is obtained from $f$ by cutting or gluing as defined above. Notice that each degree $n$ function (reduced diagram) has $2n-1$ neighbors, $n$ obtained by cutting and $n-1$ by gluing. For example if $(T,T_1,T_2)$ is a diagram of a degree $2$ function, assuming that none of the $T_i$'s are empty, its three neighbors are $(T,T_1^l,T_1^r,T_2)$, $(T,T_1,T_2^l,T_2^r)$ and $(T, \wedge\oplus^{\{1,2\}}(T_1,T_2))$.

\begin{figure}[!h]
\label{rolex}
\begin{center}
\includegraphics{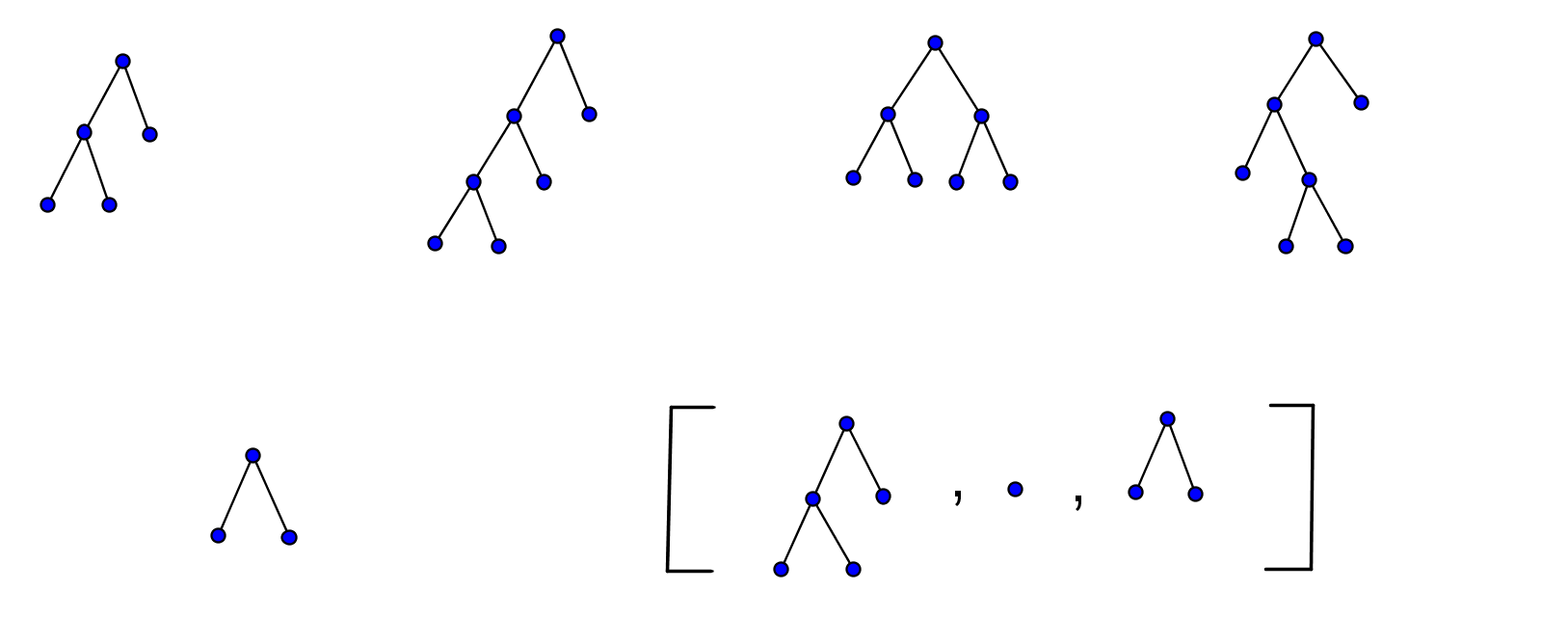}
\caption{We adopt the convention that when we write a single binary tree to denote a Thompson-like function, it means that all the binary trees on the right (that is, all except the first one) are of size $0$ (just the root). On the first line we have $L_2$ followed by its neighbors obtained by cutting. On the second one we have the two neighbors obtained by gluing.}
\end{center}
\end{figure}

\begin{definitie}
Let $\Delta=(T,T_1, \ldots, T_n)$ be an $n$-diagram and let $J \subseteq [n]$. We denote $\Delta(J)$ the diagram obtained from $\Delta$ by performing a simultaneous cutting (like in the previous definition) for any $j \in J$; more precisely for any $j \in J$ we replace $T_j$ with $(T_j^l,T_j^r)$ if $T_j \neq *$ and with $(*,*)$ if $T_j=*$, in this case also adding a caret to $T$ at the position of the unique leaf of $T_j$. The positions of the trees in the right part of the diagram are shifted naturally to right, the new diagram having degree $n+|J|$. Again, the operation does not depend on the equivalence between diagrams and naturally extends to Thompson-like functions. We denote by $g[J]$ the new function obtained from a function $g$ using this procedure. Again, $g[J]=s_{n,J} \circ g$, where $s_{n,J}$ is the unique piecewise linear bijection from $[0,n]$ to $[0,n+|J|]$, with all the slopes $1$ except on the intervals $[j-1,j]$ with $j \in J$, where the slope is $2$.
\end{definitie}

\begin{exemplu}
If $\Delta=(T,T_1,T_2,T_3)$ with all $T_i \neq *$ and $J=\{1,3\}$, then $$\Delta(J)=(T,T_1^l,T_1^r,T_2,T_3^l,T_3^r).$$
If $\Delta=(L_6,R_3,\wedge, \wedge)$ and $J=\{1,3\}$, then $$\Delta(J)=(L_6, *, R_2, \wedge,*,*).$$
If $\Delta=(T_2,*,\wedge,*)$ and $J=\{1,2\}$, then $$\Delta(J)=(T_2 \oplus^1, *, *, *,*,*).$$
\end{exemplu}

\begin{figure}[h]
\label{rolex}
\begin{center}
\includegraphics{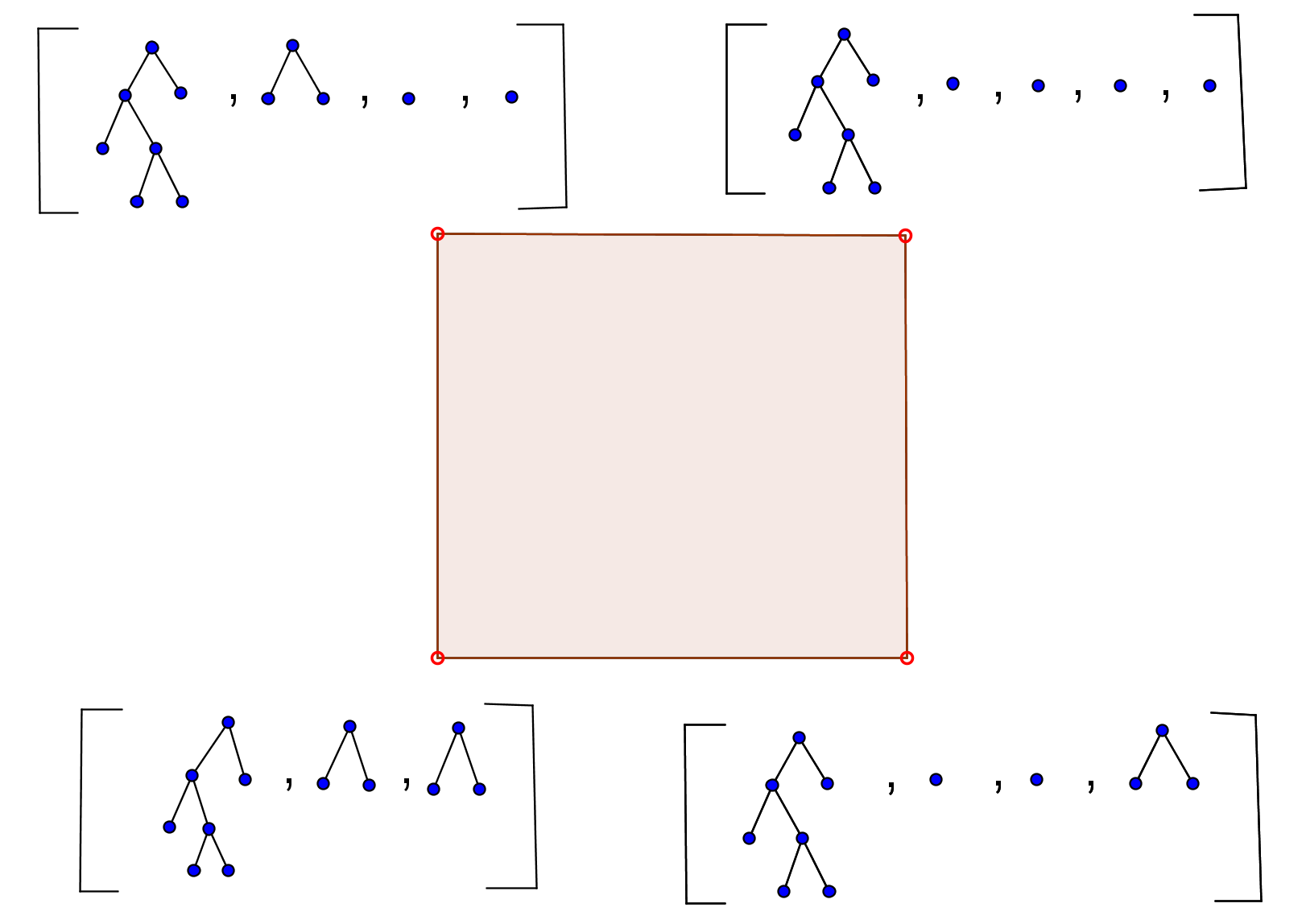}
\caption{The two dimensional cube $C(L_2 \oplus^2, \wedge, \wedge)$. }
\end{center}
\end{figure}

\begin{definitie}
We define now the set of maximal cubes which defines the cubical structure. Each cube will be labeled by an $n$-degree Thompson-like function $g$ (or, equivalently, by a reduced diagram). $C(g)$ is defined to be an $n$-dimensional cube with the vertices labeled in the set $\{ g[J] | J \subseteq [n] \}$. More precisely, the labeling is giving by the map $\phi_g: \{0,1\}^n \rightarrow F_{\infty}$, $\phi_g(\chi_n(J))=g[J]$ for any $J \subseteq [n]$, where $\chi_n(J) \in \{0,1\}^n$ is the characteristic function of $J$ in $[n]$. We denote by $X$ the cubical complex obtained by gluing the cubes along the faces with the same labels. Often, we will denote the cubes with $C(\Delta)$, where $\Delta$ is a diagram (associated with a Thompson-like function). Sometimes, we will also abuse notations and understand by $C(g)$ or $C(\Delta)$ the corresponding set of labels in $F_{\infty}$, thus we will see $F_{\infty} \subseteq X$, the set of vertices in $X$.
\end{definitie}

\begin{exemplu}
If $\Delta=(T,S) \in F$, then $C(\Delta)= \{(T,S),(T,S^l,S^r)\}$, with $(T,S)$ having coordinate $0$ and $(T,S^l,S^r)$ having coordinate $1$. If $\Delta=(T,T_1,T_2)$ with $T_1,T_2 \neq *$, then
$$ C(\Delta)=\{(T,T_1,T_2),(T,T_1^l,T_1^r,T_2),(T,T_1,T_2^l,T_2^r),(T,T_1^l,T_1^r,T_2^l,T_2^r)\}$$
with coordinates (in the order of writing): $(0,0),(1,0),(0,1),(1,1)$.
\end{exemplu}

The proof of the following Lemma is done in \cite{farley03,farley05} and also in \cite{stein} without checking the link condition (which we will do anyway). We won't repeat the proof here. If we turn upside-down the right sides of our diagrams and glue them at the bottom of the left sides, at each leaf with the same index, (after attaching transistors) we obtain Farley's complex associated with Thompson group $F$ .

\begin{lema}
$X$ is a $CAT(0)$ cubical complex.
\end{lema}

\begin{conventie}
\label{faces}
We only defined so far the maximal cubes, the ones defining the cubical structure. Of course, any face of such a cube will be a cube in its own. If $C(g)$ is a maximal cube of dimension $n$ and $I_1 \subseteq I_2 \subseteq [n]$, the following set of labels will define a face: $C(g | I_1,I_2)=\{ g[J] | I_1 \subseteq J \subseteq I_2 \}$ and all the faces of $C(g)$ are obtained in this way. If $I_1=\emptyset$, then we simply  write $C(g|I_2)$. If $g$ is represented by a diagram $(T,T_1, \ldots, T_n)$, we will write $C(T,T_1, \ldots, T_n|I_2)$.
\end{conventie}

\begin{conventie}
\label{diagrame ireductibile}
We adopt a convention for writing interior points in cubes. If $\Delta \in F_{\infty}$ is an $n$-degree diagram $(T,T_1, \ldots, T_n)$ and $x \in C(\Delta)$ is the point of coordinates $(t_1, \ldots, t_n)$ (since the cubes have an origin, there are natural coordinates), we simply represent $x$ by writing:
$$ (T,T_1, \ldots, T_n)[t_1, \ldots, t_n]$$
and call such a representation a \emph{generalized diagram} or sometimes simply, a $\emph{diagram}$.
Notice that a point $x \in X$ may have several different representations by generalized diagrams, even if we fix the initial diagram $\Delta$ to be reduced. However there is a unique generalized diagram representation with $\Delta$ reduced and all the coordinates $t_i < 1$. Indeed, we can replace the above diagram with $\Delta(J)$, where $J \subseteq [n]$ is the set of all positions where the coordinates are $1$, replacing all the coordinates $1$ with coordinates $(0,0)$ and all the others coordinates being unchanged but shifted to the right. For example $(T,T_1,T_2,T_3,T_4)[1/2,1,1/2,1/2]$ with all the $T_i \neq *$, will represent in $X$ the same element like $(T,T_1,T_2^l,T_2^r,T_3,T_4)[1/2,0,0,1/2,1/2]$. The diagram of $x$ satisfying these properties will be called \emph{the reduced (generalized) diagram} of $x$.
\end{conventie}

\begin{remarca}
The Thompson group maps cubes to cubes of same dimension, respecting the faces, thus it acts (from the left) on $X$ by cubical automorphisms (recall that we decided to respect the order of the diagram multiplication). The action is proper, but not cocompact.
\end{remarca}

\begin{conventie}
For a reduced diagram $\Delta$ we denote by $|\Delta|$ the total number of carets in $\Delta$. We call this number \emph{the norm of the diagram $\Delta$}. We also define \emph{the left norm} and \emph{the right norm}, denoted by $|.|_l$ and $|.|_r$, to be the total number of carets on the left side, respectively the right side, of the diagram. Recall that for the Thompson group $F$, Burillo defined a similar norm equal with the number of carets in a tree of the reduced diagram and that the induced distance is quasi-isometric with the Cayley distance on $F$. With our notations and identifications, the median (combinatorial) distance from an element of $F$ to the origin of the complex (which will always be the vertex representing the identity of $[0,1]$, i.e. the diagram $(*,*)$) is exactly twice the Burillo norm. As a consequence, a Cayley graph of $F$ is (quasi-)isometrically embedded in $X$ with the median distance.
\end{conventie}

\begin{figure}[h]
\label{rolex}
\begin{center}
\includegraphics{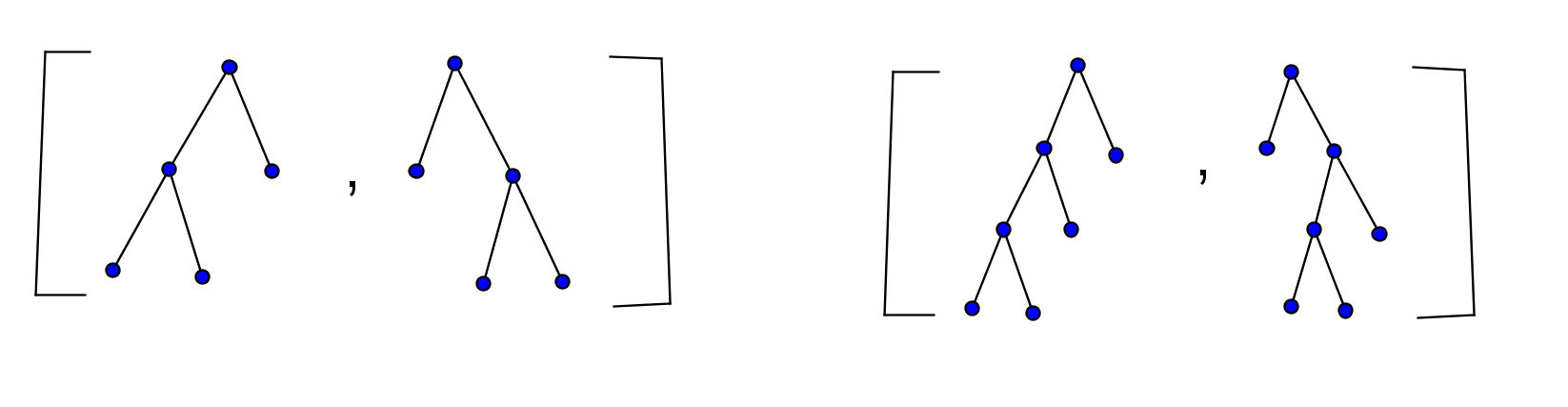}
\caption{$(L_2,R_2)$ is the first standard generator of $F$ and $(L_3,R_2 \oplus^2)$ is the element of $F$ featured in Section 2.2.}
\end{center}
\end{figure}

\subsection{Remarks}

In this section, we record some basic general facts about the complex. In the first lemma we just fix the notation for the link of the vertex. Let $v \in X$ be a vertex, represented by a $n$-diagram $\Delta$. For any $j=1, \ldots, n$, we denote by $e_j$ the edge $[\Delta,\Delta(j)]$ obtained by cutting the diagram $\Delta$ at the position $j$, and for any $k=1, \ldots, n-1$ we denote by $e_{k,k+1}$ the edge $[\Delta,\Delta_{k,k+1}]$, where $\Delta_{k,k+1}$ is the diagram obtained by gluing the diagram $\Delta$ at the positions $(k,k+1)$. The vertices of $Lk(v)$ are all the edges $e_j$, $e_{k,k+1}$.

\begin{lema}
Let $v \in X$ be a vertex. A set of edges containing $v$ belongs to a simplex in $Lk(v)$ if and only if their subscripts  are mutually disjoint as subsets of $[n]$.
\end{lema}

\noindent \emph{Proof.} By the definition of the complex, it is straightforward to check that any of the pairs $(e_j,e_{j-1,j})$ , $(e_j,e_{j,j+1})$, $(e_{j-1,j},e_{j,j+1})$, for some $j$, cannot belong to a same square. Now, if we consider a set of edges with mutually disjoint subscripts, we consider the diagram obtained from $v$ by gluing at all the positions $(k,k+1)$ for which $e_{k,k+1}$ belongs to the set. The cube of this diagram contains all the initial edges. $\clubsuit$

\vspace{2mm}

We now emphasize on several important convex subsets in $X$. Recall that $Trees$ is denoting the set of all trees (for us binary, rooted, ordered, finite trees). We can identify $Trees$ in $F_{\infty}$, namely a tree $T$ will be seen as the $|T|+1$ (reduced) diagram $(T,*, \ldots,*)$. Denote by $\mathcal{T}$ the sub-complex of $X$ consisting of all the cubes having all the vertices in $Trees$. Notice that under our identification, the cutting operation in $\mathcal{T}$ means adding a caret at the specified leaf  and the gluing operation (which is allowed only at positions labeling a caret- in order to remain in $\mathcal{T}$) means removing a caret at the specified positions, if possible. So a vertex $T \in \mathcal{T}$ has (when viewed in  $\mathcal{T}$) $|T|+1$ neighbors obtained by cutting (i.e. adding carets) and $||T||$ neighbors obtained by gluing (i.e. removing carets), where $||T||$ denotes the number of free carets in $T$. We will often write $T$ instead of $(T,*, \ldots, *)$ even when we work in $X$. The maximal cubes of $X$ which are labeled in $Trees$ define the complex $\mathcal{T}$. We have $C(T)=\{T \oplus^I | I \subseteq [n]\}$, where $n=|T|+1$. For an interior point in $C(T)$ we will simply write $T[t_1, \ldots, t_n]$ and again a point $x \in \mathcal{T}$ can be written uniquely with all the coordinates strictly less than $1$, for some tree $T$. Sometimes we will call a point in $\mathcal{T}$ a \emph{generalized tree}.

\begin{lema}
The sub-complex $\mathcal{T}$ is convex in $X$, and hence a $CAT(0)$ cubical complex in its own.
\end{lema}

\vspace{2mm}
\noindent \emph{Proof.} Since $\mathcal{T}$ is clearly connected (there is always a combinatorial path passing through the origin), it is enough to check it is locally full. Let $v \in \mathcal{T}$ be a vertex represented by a tree $T$ of size $n-1$. With the previous notations the vertices of $Lk_{\mathcal{T}}(v)$ are the edges $e_j$ for $j=1, \ldots, n$ and $e_{k,k+1}$ for any pair $(k,k+1)$ which labels a caret of $T$. It is easy to check (and follows from the analysis of the link in $X$) that a set of such vertices belongs to a simplex of $Lk_{\mathcal{T}}(v)$ if and only if their subscripts are mutually disjoint, as subsets of $[n]$. It follows immediately that $Lk_{\mathcal{T}}(v)$ is full in $Lk_X(v)$. $\clubsuit$
\vspace{2mm}

$\mathcal{T}$ is a nice model for the geometry of binary trees or, equivalently, for the geometry of dyadic partitions of a unit interval. It will be studied in details later.

\vspace{2mm}

\begin{definitie}
A tree $T \in Trees$ with $||T||=1$ (i.e. with only one free caret) is called a \emph{snake}. An infinite tree $\sigma$ is called a \emph{long snake} if it has exactly one caret at each level of depth, more precisely $|\sigma^j|=1$ (see \ref{def tree}). In other words a long snake is a (infinite) tree consisting of only mixed carets.
\end{definitie}

\begin{definitie}
For any, possibly infinite, tree $\tau$ and any non-negative integer $k$, we define $\tau[k]=\tau \cap T_k$ to be \emph{the truncation of $\tau$ at level $k$} (recall that $T_k$ is the full (finite) tree of depth $k$). Notice that if $\sigma$ is a long snake, then $\sigma[k]$ is a snake, for any $k$.
\end{definitie}

Perhaps the most ubiquitous snakes in this thesis are the most left one of depth $n$, denoted $L_n$ and the most right one of depth $n$, denoted $R_n$. $L_{\infty}=\bigcup_{n \geq 1} L_n$ is the most left long snake and $R_{\infty}=\bigcup_{n \geq 1}R_n$ is the most right one. Another useful definition is the following.

\begin{figure}[!h]
\label{rolex}
\begin{center}
\includegraphics{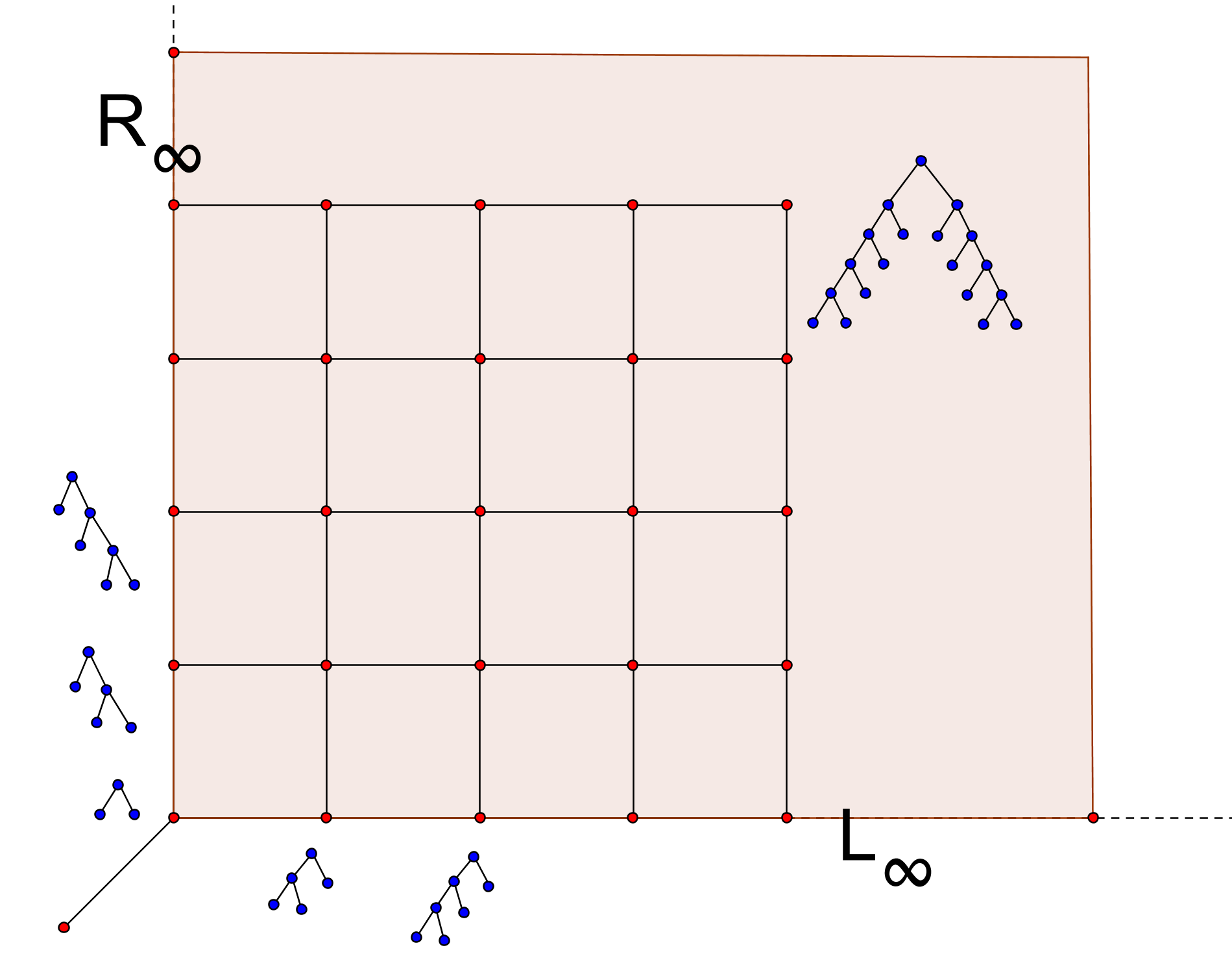}
\caption{The quadrant with tail, generated by the long snakes $L_{\infty}$ and $R_{\infty}$ (see Section 3.1 for details)}
\end{center}
\end{figure}

\begin{definitie}
\label{wings}
If $T$ is a tree, then we call the \emph{left wing} of $T$ (respectively the \emph{right wing} of $T$) to be the intersection of $T$ with $L_{\infty}$ (resp. the intersection with $R_{\infty}$). We denote by $lw(T)=|T \cap L_{\infty}|$ the size of its left wing (resp. $rw(T)=|T \cap R_{\infty}|$).
\end{definitie}

For any long snake $\sigma$ we define $c_{\sigma}: [0,\infty) \rightarrow X$, by setting $c(n)=\sigma[n]$ for any integer $n \geq 0$ and by extending the definition on each $(n,n+1)$ with the geodesic between $\sigma[n]$ and $\sigma[n+1]$. Similarly for any snake $S$ of size $m$ we define $c_S:[0, m] \rightarrow X$, by $c(k)=S[k]$ for any $k=0,\ldots,m$ (and extending by geodesics on each $[k,k+1]$).

\vspace{2mm}

\begin{lema}
\label{serpi}
For any long snake $\sigma$ and any snake $S$, we have that $c_{\sigma}$ and $c_S$ are geodesics starting at the origin $O=(*,*)$.
\end{lema}

\noindent \emph{Proof.} By truncation it is enough to work with long snakes only. Consider the one dimensional sub-complex $X_{\sigma}=\bigcup_{n=0}^{\infty}[\sigma[n],\sigma[n+1]]$. Since it is clearly connected and for any $n \geq 0$, $Lk_{X_{\sigma}}(\sigma[n])$ consists of only two non-linked vertices, $X_{\sigma}$ is convex and we have the claim. $\clubsuit$

\vspace{2mm}

We now give a formula for the distance in the complex $\mathcal{T}$ which reduces a bit the complexity of computation. Recall that the $1$-skeleton of a product of two cubical complexes $Y$ and $Z$ is just the graph cartesian product of the $1$-skeletons of the two complexes and the $CAT(0)$ distance is given by $\sqrt{d_Y^2+d_Z^2}$. We denote by $\mathcal{T}^n$ the $n$-times direct product of the complex $\mathcal{T}$. The maximal cubes in $\mathcal{T}^n$ are labeled by the $n$-tuples of trees. More precisely, the vertices of the cube labeled by $(T_1, \ldots, T_n)$ are:
$$\{(T_1\oplus^{I_1}, \ldots, T_n \oplus^{I_n})| I_1 \subseteq [|T_1+1], \ldots, I_n \subseteq[|T_n|+1]\}$$

For any tree $T$ of size $m-1$ we define the injective cubical map $\phi_T:\mathcal{T}^m \rightarrow \mathcal{T}$, called \emph{lantern}, defined on the vertices by $$\phi_T((A_1, \ldots, A_m))=T \oplus (A_1, \ldots, A_m)$$

Notice that the notations from $\ref{conventie lipire}$ make sense also if we attach generalized trees, so if $A_1, \ldots, A_m \in \mathcal{T}$ are generalized trees, then above we have the full definition of the lantern $\phi_T$.
We also have an obvious inclusion relation for (generalized) trees, so $\phi_T$ is in fact a bijective cubical map from $\mathcal{T}^m$ to the space $\mathcal{T}_T=\{x \in \mathcal{T} | T \subseteq x\}$.

\begin{lema}
\label{calcul recursiv}
(Recursive computing in Trees) Let $T$ be a tree of size $m-1$ and $A_1, \ldots, A_m$ and $B_1, \ldots, B_m$ be some (generalized) trees. Then we have that $\phi_T$ is distance preserving and in particular:

$$ d(T \oplus(A_1, \ldots, A_m), T \oplus (B_1, \ldots, B_m))= \sqrt{\sum_{k=1}^m d(A_k,B_k)^2}.$$
\end{lema}

\vspace{3mm}

\noindent \emph{Proof.} Since the lantern $\phi_T$ is cubical and injective, we only have to check that its image, the sub-complex $\mathcal{T}_T$ is convex in $\mathcal{T}$. First we notice that $\mathcal{T}_T$ is connected, since, for example, between any two points there is a combinatorial path passing through $T$. It remains to check that $\mathcal{T}_T$ is locally full in $\mathcal{T}$. Let $S \supset T$ be a tree. Denote by $A$ the set of all pairs $(k,k+1)$ which index a free caret of $S$ which is also a free caret of $T$ (when $S$ and $T$ are viewed as embedded in the full binary tree, i.e. the caret under discussion has the same position in the full binary tree (\ref{def tree}) and belongs as a free caret in both $T$ and $S$). Denote by $B$ the set of all pairs $(l,l+1)$ which index all the other free carets of $S$, i.e. those who are not free carets in $T$. A simplex in $Lk_{\mathcal{T}}(S)$ with all the vertices in $Lk_{\mathcal{T}_T}(S)$ consists of some edges (in $\mathcal{T}$) of the form $e_j$ and some edges of the form $e_{k,k+1}$ with $(k,k+1) \in B$, with all the subscripts involved disjoints as subsets of $[|S|+1]$. Denote by $V$ the tree obtained by deleting from $S$ all the carets $(k,k+1) \in B$ corresponding to the edges $e_{k,k+1}$ in the simplex. The tree $V$ still contains $T$ and the cube indexed by $V$ is a cube in $\mathcal{T}_T$ and contains all the edges involved in the simplex. $\clubsuit$

\newpage

\begin{center}
$\clubsuit$
\end{center}

\newpage

\section{Isometries}

In this chapter we prove Theorem 1., Theorem 2. and Corollary 1. from the introduction. In section 2.1 we show that some elements of $F$, which can be viewed as small perturbation of the identity, are hyperbolic and we use them to show the existence of the flats. In section 2.2. we study the behavior of a particular isometry in $F$. In section 2.3 we prove the first half of Theorem 2, namely we compute the translation lengths. In section 2.4 we complete the proof by showing the second half of the theorem: we show the existence of parabolic isometries in $F$.

\subsection{Flats}

Warning: in this section we use freely the standard basic facts about hyperbolic isometries in \cite{book}.

\begin{definitie}
Let $S$ be a snake with the leaves of its only free caret labeled $k,k+1$. The element $(S \oplus^k,S \oplus^{k+1}) \in F$ is called \emph{ a small perturbation of the identity} or \emph{ a rotation} based on the snake $S$.
\end{definitie}

Notice that rotations are just scaled copies of the first standard generator of $F$, namley $(L_2,R_2)$ on some dyadic interval (specified by a snake). We choose the name rotation since in all the diagrams representing rotations, one tree is the rotation of the other tree in the sense of Sleator-Tarjan-Thurston \cite{rotation}. The famous conjecture of proving that the diameter of the associaedron in $2n-6$ for $n \geq 11$ can be stated in terms of The Thompson Group equipped with the set of generators consisting of all rotations as we defined them.

\begin{lema}
The rotations are all hyperbolic isometries with translation length $\sqrt{2}$.
\end{lema}

\vspace{2mm}
\noindent \emph{Proof.} Let $g=(S \oplus^k,S \oplus^{k+1})$ be a rotation based on a snake $S$. It is enough to find a point $x \in X$ such that $d(gx,g^{-1}x)=2d(gx,x)$. We choose $x=S \oplus^{\{k,k+1\}}$ and we have $gx=S \oplus^k L_2$ and $g^{-1}x=S \oplus^{k+1}R_2$. Computing with  the recursive  formula \ref{calcul recursiv} we easily have $d(gx,g^{-1}x)=2d(gx,x)=2\sqrt{2}$. $\clubsuit$
\vspace{2mm}

\noindent \emph{Proof of Theorem 2.} We use the following fact, which can be easily proved by induction starting with the flat strip lemma: if $f_1,f_2, \ldots, f_m$ are two by two commuting hyperbolic isometries with the same translation length and with some given axis $L_1,L_2, \ldots, L_m$ two by two orthogonal and all meeting in some point $x$ then $conv(L_1\cup L_2 \ldots \cup L_m)$ is isometric with an $m$-dimensional flat.

Let $n \geq 1$ be a natural number and let $T=T_n$ be the full tree of depth $n$. Let $f_j=(T \oplus^{2j-1},T\oplus^{2j})$ for $j=1, \ldots, 2^{n-1}$. Notice that the $f_j$'s are all rotations (although not written in the reduced form) and commute two by two (they have disjoint support). Consider $x=T \oplus^M$, where $M=\{2,4,6, \ldots, 2^n\}$.

A simple computation shows that $f_jx=T \oplus^{(M-\{2j\})\cup \{2j-1\}}$, for any $j=1, \ldots, 2^{n-1}$.  Notice that all the elements $x, f_1x, \ldots, f_{2^{n-1}}x$ belong to the $2^n$ dimensional cube $C(T)$, hence we have $d(f_jx,x)=\sqrt{2}$ for any $j$ and $d(f_jx,f_ix)=2$ for any $i \neq j$. In particular, by Lemma 2.1.2, the sets $L_j=\bigcup_{k \in \mathbb{Z}}[f_j^{k}x,f_j^{k+1}x]$ are axes for the $f_j$'s. By Pythagoras Theorem we have $\angle_x(L_i,L_j)=\pi/2$ and we find a $2^{n-1}$ dimensional flat. $\clubsuit$

\begin{remarca}
Notice that in the proof above we can vary the tree $T$ and consider snakes based at the free carets of $T$. In this way one can obtain other variants of flats.
\end{remarca}

\begin{figure}[!h]
\label{rolex}
\begin{center}
\includegraphics{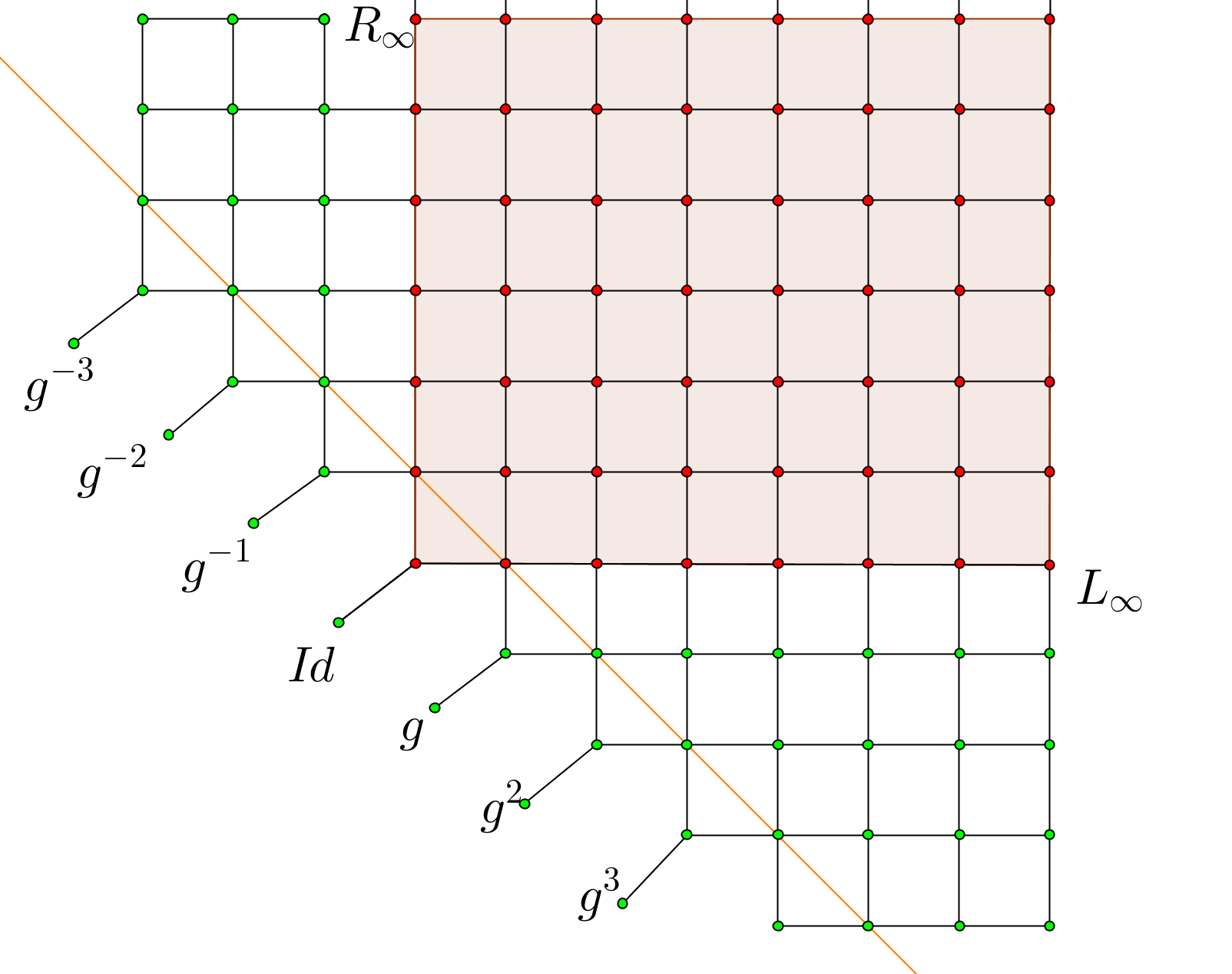}
\caption{A two dimensional sheet at the "bottom" of the complex, near the origin and around the quadrant generated by $L_{\infty}$ and $R_{\infty}$.The red vertices denotes binary trees, more precisely the quadrant generated by $L_{\infty}$ and $R_{\infty}$ in the previous figure. The green vertices represents elements near the quadrant which are not binary trees. The leaves of the figure are indexed by the cyclic group generated by the first standard generator, denoted $g$. The axis of $g$ is the orange line: it goes through the red vertices $L_2$ and $R_2$ and through the green vertices of the form $(L_n,*,*,R_{n-2})$ (positive powers) and $(R_n,L_{n-2},*,*)$ (negative powers) for $n \geq 3$.
}
\end{center}
\end{figure}

\subsection{A particular isometry}

This section illustrates the simple strategy used to prove Theorem 1, avoiding many technical issues. We study the third most simple (to draw) element of $F$, namely $g=(L_3,R_2 \oplus^2)$. We show that the translation length of $g$ is $\sqrt{5}$ and the minimal displacement on vertices is $\sqrt{6}$. We also give examples of points where the displacement of $g$ is arbitrarily close to $\sqrt{5}$. The section is independent from the rest of the text and can be skipped, although we do not recommend this.

\begin{propozitie}
The translation length of $g$ is $\sqrt{5}$.
\end{propozitie}

\noindent \emph{Proof.} A simple computation shows that $g^n=(L_{2n+1}, K_n)$, where $K_n=R_{n+1} \oplus^{[2,n+1]}$ (in this setting we denote $[i,j]=\{i,i+1, \ldots, j\}$ for two positive integers $i < j$).

\vspace{2mm}

\noindent \emph{The upper bound.}

\vspace{2mm}

By Lemma \ref{asimptotic} in the next section and also using $\ref{calcul recursiv}$ and $\ref{serpi}$ during the computation, for any $n$, we have:
\begin{eqnarray*}
n|g|=|g^n| &\leq& d(g^nK_n,K_n)\\
           &=&d(L_{2n+1},K_n) \\
           &=& \sqrt{d(O,L_{2n})^2+d(O,R_n\oplus^{[1,n]})^2}\\
           & \leq & \sqrt{4n^2+(d(O,R_n)+d(R_n,R_n \oplus^{[1,n]})^2}\\
           &=& \sqrt{4n^2 + (n + \sqrt{n})^2}\\
           &=& \sqrt{5n^2+ 2n^{3/2}+ n}.
\end{eqnarray*}

\noindent When $n \rightarrow \infty$, we have $|g| \leq \sqrt{5}$.

\vspace{2mm}

\noindent \emph{The lower bound.}

\vspace{2mm}

We implement the asymptotic formula from Lemma \ref{asimptotic}, starting to iterate at the origin. We have thus to estimate $d(0,(L_{2n+1},K_n))$. By the Lemma \ref{proiectie trees} in the next section, the projection of $(L_{2n+1},K_n)$ to $\mathcal{T}$ is $L_{2n+1}$ and then the triangle $O,L_{2n+1},(L_{2n+1},K_n)$ is obtuse at $L_{2n+1}$. By also using the fact that the projection of $K_n$ on $R_{\infty}$ is $R_{n+1}$ (Lemma \ref{proiectie serpi} in the next section), we have:
\begin{eqnarray*}
d(g^n(O),O)^2 &\geq& d(O,L_{2n+1})^2 + d((L_{2n+1},K_n),L_{2n+1})^2 \\
              &=& d(O,L_{2n+1})^2 + d(g^nO,g^nK_n)^2 \\
              &=& d(O,L_{2n+1})^2 + d(O,K_n)^2 \\
              &\geq& (2n+1)^2 + d(\Pi_{R_{\infty}}(0),\Pi_{R_{\infty}}(K_n))^2 \\
              &=& (2n+1)^2 + d(0,R_{n+1})^2 \\
              &=& 5n^2 + 6n + 2.
\end{eqnarray*}

By the asymptotic formula we have $|g| \geq \sqrt{5}$.  $\clubsuit$

\vspace{2mm}

We now exhibit some points where the displacement is arbitrary close to $\sqrt{5}$. For any $n$ we consider $x_n=(R_{n+2} \oplus^2)(0,0,0, \frac{n-1}{n},\frac{n-2}{n}, \ldots, \frac{1}{n}, 0,0)$. A direct computation with \ref{calcul recursiv} shows that $$d(gx_n,x_n)=\sqrt{4+ 1/n^2+ 1/n^2+ \ldots + 1/n^2+1}=\sqrt{5+1/n}.$$
This, of course, gives another proof that $|g| \leq \sqrt{5}$. By Theorem 1, $\sqrt{5}$ is not achieved and we will see that , indeed, $g$ is parabolic. A skeleton version of the argument is given by the proof of the next Proposition.

\vspace{2mm}

\begin{propozitie}
For the isometry $g$, the minimal displacement on vertices is $\sqrt{6}$.
\end{propozitie}

\vspace{2mm}

\noindent \emph{Proof.} First of all let us notice that if $w=R_3 \oplus^2$ then a straightforward computation shows that  $gw=L_3 \oplus^4$. By applying $\ref{calcul recursiv}$ with the lantern $\phi_{R_2}$ we have that

$$ d(gw,w)= \sqrt{d(O,L_2)^2 + d(O,\wedge)^2 + d(O, \wedge)^2}=\sqrt{6}$$

We are left to show that for any vertex $v \in X$ we have $d(gv,v) \geq \sqrt{6}$.

\vspace{2mm}

Let $v=(A,A_1, \ldots, A_n)$ be a degree $n$ vertex in $X$. We assume for the moment $n \geq 2$ (the case $n=1$ is easier but a bit different so we postpone it until the end of this proof). We reduce the computation $d(gv,v)$ by multiplying both terms with a convenient isometry $ h \in F$, that is estimating instead the (same) distance $d(hgv,hv)$. We define $h=(L_{n-1} \oplus (A_1, \ldots, A_n),A)$. We clearly have $hv=L_{n-1}$. We cannot determine the precise form of $hgv$ but the following easy to notice information will be enough: the first slope of $hgv$ interpreted as piecewise linear map is $2^{n+1}$ and the last slope is $1$ (by looking at the derivatives of the functions at $0$ and $1$). If $(T,T_1, \ldots, T_n)$ is the reduced diagram of $hgv$, the information on the slopes reads: $lw(T)=lw(T_1)+n+1$ and $rw(T)=rw(T_n)$ (recall \ref{wings}). So it is left to show that the distance between a diagram with these two properties and $L_{n-1}$ is at least $\sqrt{6}$.

\vspace{2mm}

By applying the Lemma \ref{proiectie trees} in the next section we see that the triangle $(T, T_1, \ldots, T_n)\\,T,L_{n-1}$ is obtuse at $T$. Denoting by $k=(L_{n-1}\oplus (T_1, \ldots, T_n),T) \in F$ and applying \ref{calcul recursiv}  we have the following estimation:

\begin{eqnarray*}
d((T,T_1, \ldots, T_n),L_{n-1})^2 &\geq& d(T,L_{n-1})^2 + d((T,T_1, \ldots,T_n),T)^2 \\
              &=& d(T,L_{n-1})^2 + d(k(T,T_1, \ldots,T_n),kT)^2 \\
              &=& d(T,L_{n-1})^2+d(L_{n-1},L_{n-1}\oplus (T_1, \ldots, T_n))^2\\
              &=& d(T,L_{n-1})^2 + d(O,T_1)^2+ \ldots +d(O,T_n)^2
\end{eqnarray*}

Let us analyze the situation. If $T_1 \neq *$ then $d(T, L_{n-1}) \geq d(\Pi_{L_{\infty}}(T),L_{n-1}) \geq 3$, since $lw(T) \geq n+2$ in this case. So we may assume that $T_1=*$. If $rw(T_n) \geq 2$ then $d(O,T_n) \geq d(\Pi_{R_{\infty}}(O),\Pi_{R_{\infty}}(T_n)) \geq 2$, so our main inequality becomes:

\begin{eqnarray*}
d((T,T_1, \ldots, T_n),L_{n-1})^2 &\geq& d(T,L_{n-1})^2 + d(O,T_n)^2 \\
                                  &\geq& d(\Pi_{L_{\infty}}(T),\Pi_{L_{\infty}}(L_{n-1}))^2 + 4 \\
                                  &\geq& 4 + 4 = 8
\end{eqnarray*}

So we may assume that $T_1=*$ and $T_n=\wedge$, that is, $lw(T)=n+1$ and $rw(T)=1$. By the Caret Counting Lemma \ref{numarare} in the next section, we have $$n+1 \leq |T|= n-1 + |T_1| + |T_2| + \ldots + |T_n|$$

From this equation at least one of the trees $T_2, \ldots, T_{n-1}$ must be non-empty (this case is possible only if $n \geq 3$, of course). If $T_k$ denotes the non-empty tree, our main inequality shows:
\begin{eqnarray*}
d((T,T_1, \ldots, T_n),L_{n-1})^2 &\geq& d(T,L_{n-1})^2 + d(O,T_n)^2 + d(O,T_k)^2\\
                                  &\geq& 4 + 1 +1 = 6
                                  \end{eqnarray*}

To finish we only have to deal with the case $n=1$, that is when $v \in F$. We have $d(gv,v)=d(v^{-1}gv,O)$, so we only have to estimate the distance from $(T,S)$ to the origin where $(T,S)$ is an reduced diagram with $lw(T)=lw(S)+2 \geq 3$ and $rw(S)=rw(T)+1 \geq 2$. With the same observations as before, we have:

\begin{eqnarray*}
d((T,S),O)^2 &\geq& d((T,S),T)^2 + d(O,T)^2 \\
                                  &=& d(S,O)^2 + d(O,T)^2 \\
                                  &\geq& rw(S)^2 + lw(T)^2 \\
                                  &\geq& 4 + 9 =13.
\end{eqnarray*}
$\clubsuit$
\vspace{2mm}

\subsection{Translation Lengths}

In this section we show the first half of Theorem 1, namely:

\begin{propozitie}
\label{lungime translatie}
If $ g \in F$ is an irreducible element (that is, $g$ only fixes $0$ and $1$ as a homeomorphism of the unit interval), then its translation length in the model is $\sqrt{log_2(g'(0))^2 + log_2(g'(1))^2}$.
\end{propozitie}

We start with a string of some quite general lemmas. The first one is an unpublished remark of N. Monod.

\begin{lema}
\label{asimptotic}
\emph{(The asymptotic formula)} If $M$ is a  $CAT(0)$ space and $f$ is an isometry of $M$, then for any $x \in M$ we have $|g|=lim_{n \rightarrow \infty}\frac{d(g^n(x),x)}{n}$. In particular, for any $n \in \mathbb{N}$ we have $|g^n|=n|g|$.
\end{lema}

\noindent \emph{Proof.} It follows from the triangle inequality that the sequence $d(g^nx,x)$ is subadditive and hence the limit exists. The triangle inequality further implies that the limit is independent of $x$ and is bounded above by $|g|$.

For the reverse inequality, it is enough to prove by induction on $n$ that $d(g^{2^n}x,x) \geq 2^n|g|$ holds for all $x \in M$. Consider the midpoint $x'$ of $[x,g^{2^{n-1}}x]$. Since $g^{2^{n-1}}x'$ is the midpoint of $[g^{2^{n-1}}x,g^{2^n}x]$, the $CAT(0)$ inequality implies $d(g^{2^n}x,x) \geq 2d(g^{2^{n-1}}x',x')$. We conclude by induction since the inequality holds for $n=0$ by the definition of $|g|$. $\clubsuit$

\begin{lema}
\label{interplay}
\emph{(The interplay between CAT(0) and median)} Let $M$ be a CAT(0) cubical complex. For any two vertices $x,y \in M$ consider the convex sub-complex $Conv[x,y]$ consisting of all the cubes with all the vertices in the median (i.e. combinatorial) segment between $x$ and $y$. The image of the $CAT(0)$ geodesic from $x$ to $y$ is included in $Conv[x,y]$. Moreover, $Conv[x,y]$ can be realized as a sub-complex of the integer cube structure of $\mathbb{R}^N$, where $N$ is the maximal dimension of a cube in $Conv[x,y]$, more precisely: $Conv[x,y]$ is cubical isomorphic (hence median and $CAT(0)$ isometric) with a (in general non-convex) sub-complex of $\mathbb{R}^N$, endowed with the shortest path (inside the sub-complex) metrics. Similarly, one can define $Conv[C_1,C_2]$ for two cubes $C_1$ and $C_2$ to be the sub-complex consisting of all cubes with all vertices in a median segment between two vertices, one in $C_1$ and the other one in $C_2$. $Conv[C_1,C_2]=Conv[x,y]$ for some $x \in C_1$ and $y \in C_2$.
\end{lema}

\noindent \emph{Proof.} \cite{ardila} Proposition 3.2. (see also the introduction) $\clubsuit$

\vspace{3mm}

\begin{remarca}
We won't need the previous lemma at full power. The bound $N$ is not important for us and the considerations about the metrics also.
\end{remarca}

\begin{lema}
\label{geodezice separate vf}
\emph{(Separated geodesics near a vertex)} Let $M$ be a CAT(0) cubical complex and let $C_1$ and $C_2$ be two cubes such that their intersection consists of exactly one vertex $V$. Let $c_1$ and $c_2$ be two geodesics starting at $V$ and some $\epsilon >0$ such that $c_1[0,\epsilon] \subseteq C_1$ and $c_2[0,\epsilon] \subseteq C_2$. Then we have that $ \angle_V(c_1,c_2) \geq \pi/2$.
\end{lema}

\noindent \emph{Proof.} We realize $Conv[C_1,C_2]$ in $\mathbb{R^N}$ for some $N$, like in the Lemma \ref{interplay}. We may assume that $V$ is mapped into the origin of $\mathbb{R}^N$. For any $ 0 \leq t \leq \epsilon$, let $(\alpha_1, \ldots, \alpha_N)$ be the coordinates of $c(t)$ and $(\beta_1, \ldots, \beta_N)$ the coordinates of $c'(t)$ (after identification). By the assumption that the two cubes intersect only at $V$, we have that for any $ 1 \leq i \leq N$, if $\alpha_i$ and $\beta_i$ are both non-zero then they have opposite signs. It follows that:

$$d(c(t),c'(t))^2 \geq d_N(c(t),c'(t)^2= \sum_{i=1}^n(\alpha_i-\beta_i)^2 \geq \sum_{i=1}^n \alpha_i^2+ \sum_{i=1}^n \beta_i^2 = 2t^2$$

\noindent where $d_N$ denotes the euclidian distance on $\mathbb{R}^N$. Letting $t \rightarrow 0$ we get the conclusion. $\clubsuit$

\vspace{2mm}

\begin{definitie}
\label{terminatie}
A sequence of trees is a list of trees $\Delta=(T_1, \ldots, T_n)$ whose leaves are indexed from left to right from $1$ to $|T_1|+ \ldots + |T_n|+n$ (for example, the leaf numbered $|T_1|+2$ will be the most left leaf of $T_2$ etc.) The \emph{termination} of $(T_1, \ldots, T_n)$ is the unique sequence of trees $\Delta'=(S_1, \ldots, S_m)$ with the same total number of leaves, with each $S_j=*$ or $S_j=\wedge$ and such that the leaves indexed by  $k$ and $k+1$ belongs to a caret in $\Delta$ if and only if the leaves $k$ and $k+1$ belongs to a caret in $\Delta'$. So basically the termination of a sequence of trees just remember which leaves belongs to a same caret and which not.
\end{definitie}

\begin{lema}
\label{proiectie trees}
\emph{(Vertex projection on Trees)} Let $V \in X$ be a vertex with reduced diagram $(T,T_1, \ldots, T_n)$. Then the projection of $V$ on the subcomplex $\mathcal{T}$ is $T$.
\end{lema}

\noindent \emph{Proof.} We show that $\angle_T(c,c') \geq \pi/2$, where $c$ is the geodesic from $T$ to $V$ and  $c'$ is any starting at $T$ with the image in $\mathcal{T}$. Taking any sufficiently small $t>0$, we shall determine to which cubes $c(t)$ and $c'(t)$ belongs.
Lemma \ref{interplay} and a simple analysis of the median segment from $T$ to $V$  show that $c(t)$ belongs to the cube $C_1=C(T, \Delta_1|I)$ (see \ref{faces}), where $\Delta_1$ is the termination of the trees sequence $(T_1, \ldots, T_n)$ (see \ref{terminatie}) and $I$ is the set of positions in the sequence where we have a $\wedge$ in $\Delta$ (this is the cube determined by all the vertices in the median segment from $T$ to $V$, at distance one from $T$). $c'(t)$ belongs to a cube $C_2=C(T, \Delta_2|J)$, where $\Delta_2$ consists only of $*$ and $\wedge$ and the leaves of the $\wedge$ trees must be labeled with $(l,l+1)$, positions at which we have a caret in $T$ (otherwise we are not in $\mathcal{T}$ anymore). $J$ consist of the positions where we have a $\wedge$ and maybe some others (this choice of $C_2$ covers all the possibilities for the first cube visited by a geodesic in $\mathcal{T}$ starting at $T$). Notice also that $(T, \Delta_2)$ is not written necessarily in the reduced form. Nevertheless, since the diagram $(T,T_1, \ldots, T_n)$ is reduced, we have that $C_1 \cap C_2=\{T\}$. By Lemma \ref{geodezice separate vf} we have $\angle_T(c,c') \geq \pi/2$. $\clubsuit$

\begin{lema}
\label{proiectie serpi}
\emph{( Vertex Projection on Snakes)} Let $T$ be a vertex in $Trees$ and let $\sigma$ be a long snake. The projection of $T$ on the half-line geodesic associated with $\sigma$ is $T \cap \sigma$.
\end{lema}

\noindent \emph{Proof.} Let $n$ be such that $\sigma[n]= T \cap \sigma$. Let $c$ be the geodesic from $\sigma[n]$ to $T$, $c_1$ the geodesic from $\sigma[n]$ to $\sigma[n-1]$ and $c_2$ the geodesic from $\sigma[n]$ to $\sigma[n+1]$. We have to show that $\angle_{\sigma[n]}(c,c_1) \geq \pi/2$ and $\angle_{\sigma[n]}(c,c_2) \geq \pi/2$. We focus on the first angle.

By Lemma $\ref{interplay}$, the geodesic $c$ starts in the cube $C=C(\sigma[n]|I)$, where $I$ is the maximal subset of $[n+1]$ such that $ \sigma[n] \oplus^I \subseteq T$. $c_1$ starts, of course, in the cube $C_1=\sigma[n-1]\oplus^k$, where $k$ is the position at which we add a caret to obtain $\sigma[n]$. Clearly $C_1 \cap C_2= \{\sigma[n]\}$ and by Lemma \ref{geodezice separate vf}, we have the angle inequality.

The second inequality follows in the same way, noticing that the cube $C$ doesn't contain $\sigma[n+1]$. $\clubsuit$

\vspace{2mm}

We are now in the position to prove Proposition 2.3.1. Let $g$ be an irreducible element of $F$. Notice that if $g'(0) > 1$ then $g'(1) < 1$ and vice-versa. By replacing, if necessary, $g$ with $g^{-1}$ we can assume $g'(0) >1$. Denote $\alpha=|log_2(g'(0))|$ and $\beta=|log_2(g'(1))|$.

\vspace{2mm}

\noindent \emph{The lower bound.}

\vspace{2mm}

For any $n$, let $(T_n,S_n)$ be the reduced diagram of $g^n$. We have $lw(T_n)=lw(S_n)+n\alpha$  and $rw(S_n)=rw(T_n)+ n\beta$. The lower bound follows from the next lemma via \ref{asimptotic}.

\vspace{2mm}

\begin{lema}
\label{kk}
Let $(T,S)$ be a reduced diagram. Then we have:
$$ d((T,S),O) \geq \sqrt{lw(T)^2 + rw(S)^2}$$
\end{lema}

\begin{figure}[!h]
\label{rolex}
\begin{center}
\includegraphics{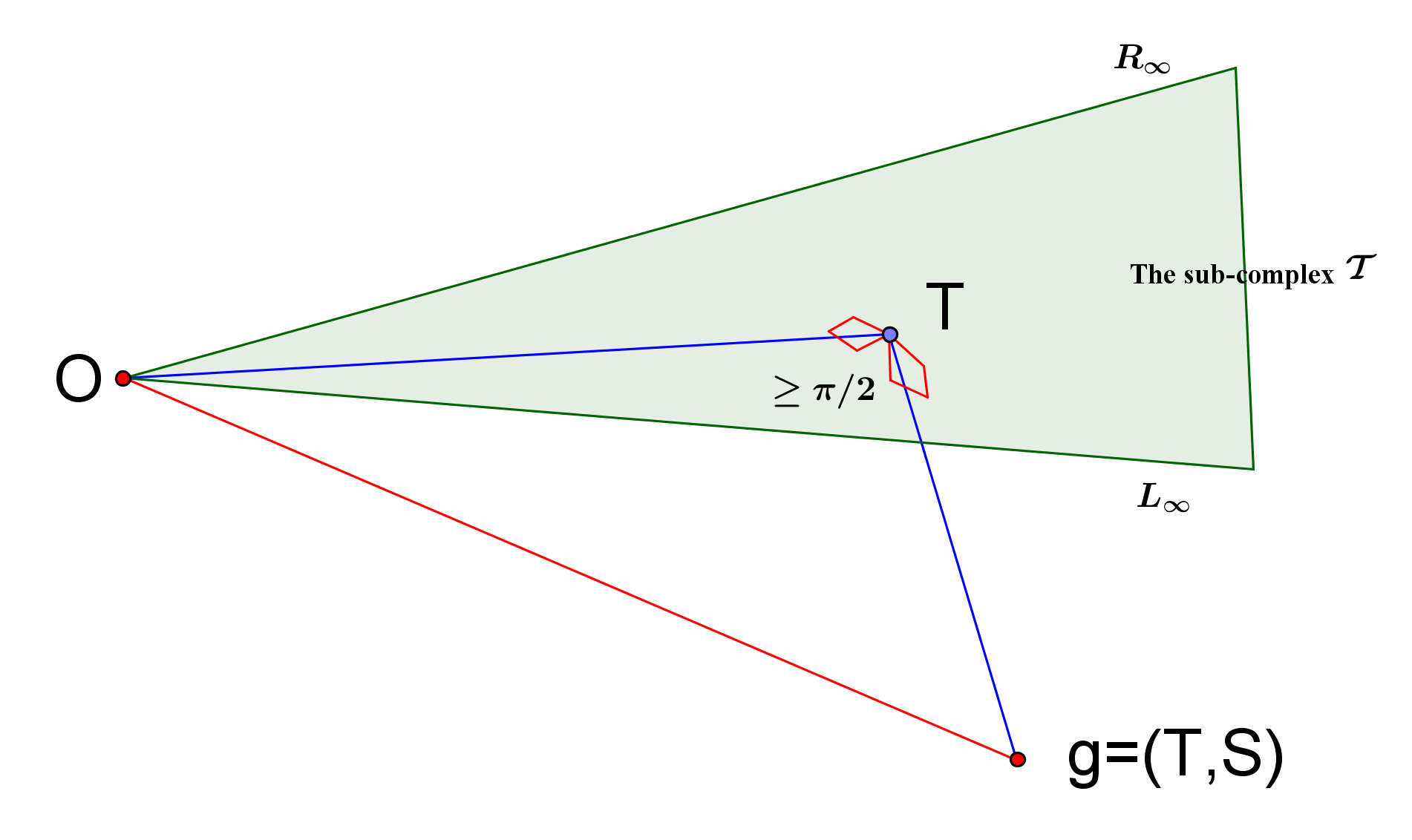}
\caption{In the proof of the lower bound (Lemma \ref{kk}), we first break the diagram in two. The geodesics $[O,T]$ and $[T,(T,S)]$ just kiss to give an obtuse triangle.}
\end{center}
\end{figure}
\noindent \emph{Proof.} The triangle $(T,S),T,O$ is obtuse at $T$ by \ref{proiectie trees}. By using also \ref{proiectie serpi}, we have:
\begin{eqnarray*}
d((T,S),O)^2 &\geq& d((T,S),T)^2 + d(T,O)^2\\
           &=& d((T,S)O,(T,S)S)^2+ d(T,O)^2 \\
           &=& d(O,S)^2 + d(T,O)^2\\
           &\geq& d(\Pi_{R_{\infty}}(O),\Pi_{R_{\infty}}(S))^2+d(\Pi_{L_{\infty}}(T),\Pi_{L_{\infty}}(O))^2\\
           &=& rw(S)^2 + lw(T)^2.
           \end{eqnarray*}
 $\clubsuit$

\begin{figure}[!h]
\label{rolex}
\begin{center}
\includegraphics[scale=0.5]{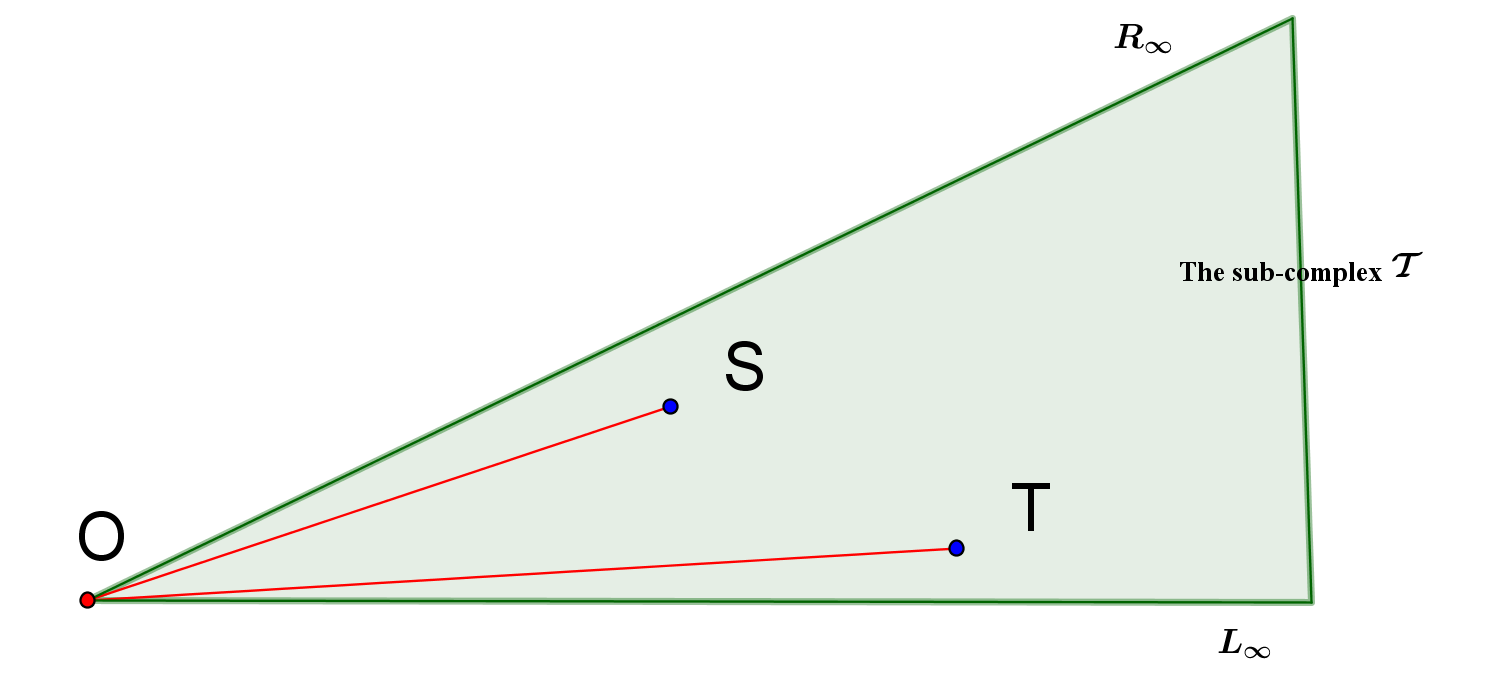}
\caption{This figure illustrates the continuation of the argument in Lemma \ref{kk} (see the previous picture). Just says that $d((T,S),T)=d(O,S)$ by applying the isometry $g$. So we know have $d(g,O)^2 \geq d(O,T)^2 + d(O,S)^2$...}
\end{center}
\end{figure}

\begin{figure}[!h]
\label{rolex}
\begin{center}
\includegraphics{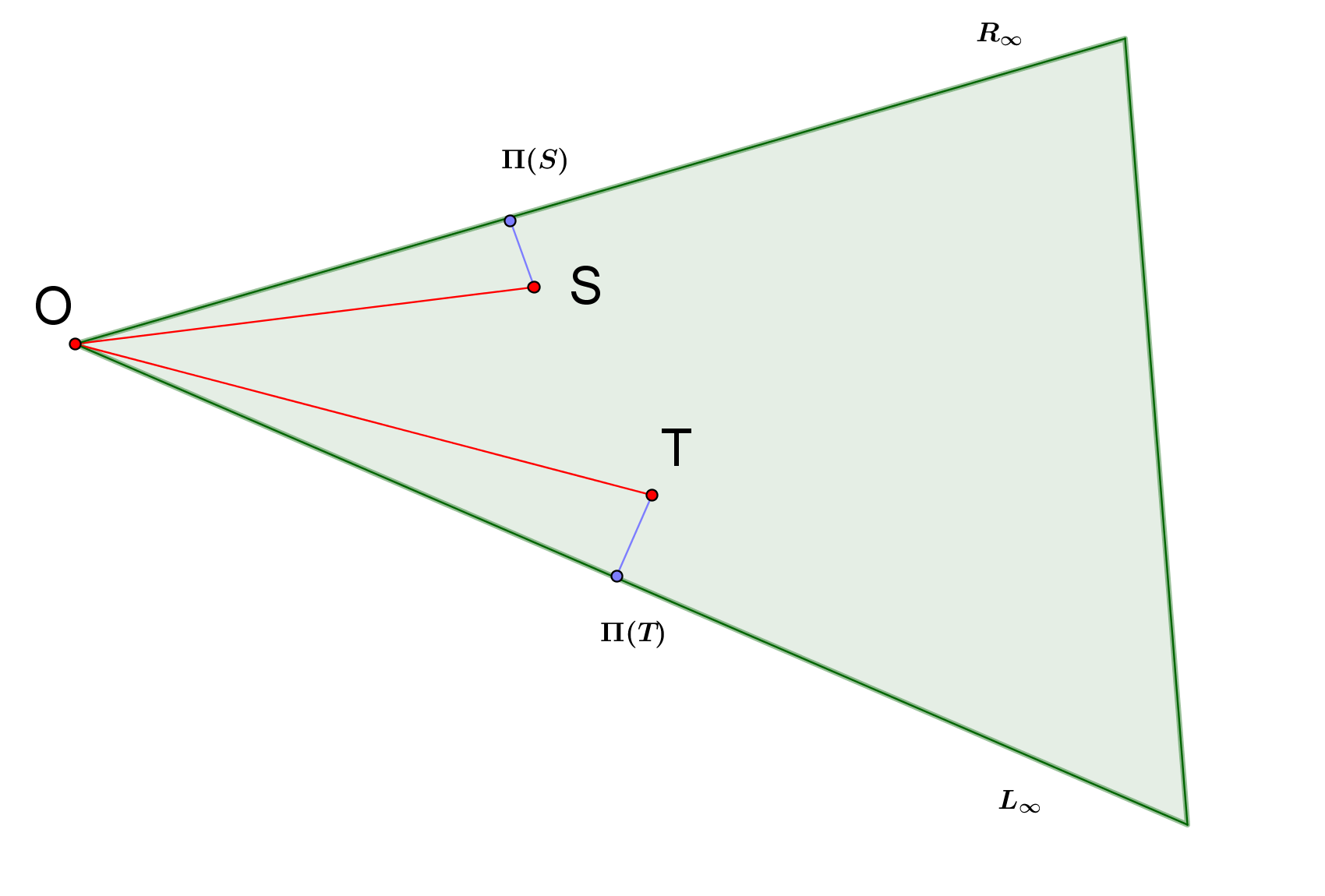}
\caption{The last step in the proof of Lemma \ref{kk}: we project $S$ on $R_{\infty}$ and  $T$ on $L_{\infty}$ to get two new obtuse triangles( we denote both projection with $\Pi()$ since there is no confusion). On $L_{\infty}$ we read the value of the first slope and on $R_{\infty}$ the value of the last slope. We get the lower bound, without using the other two edges of the triangles, which are not on $L_{\infty}$ or $R_{\infty}$. We anticipate and use this figure to also explain the parabolicity result (see Section 2.4 for details). No matter at which point $x \in X$ we evaluate $d(gx,x)$ we will reduce the computation to a configuration similar with the one in this figure. The condition $g'(0)g'(1) \neq 1$ implies combinatorially that in the configuration there is always "a missing caret", i.e. a caret not contained on $L_{\infty}$ or $R_{\infty}$. Geometrically this means that at least one of the two obtuse triangles will be non-degenerated and the inequality is strict.}
\end{center}
\end{figure}

\noindent \emph{The upper bound.}

\vspace{2mm}

In order to prove the upper bound, we need to estimate the shape of the diagrams of  the iterations $g^n$. If $(T,S)$ is a diagram of $g$, we denote $U=T \cap S$. Notice that $|U| \geq 1$.

Since it is simpler and quite illuminating we first treat the case $|U|=1$, that is $U= \wedge$. In this case $T=\wedge \oplus(A,*)$ and $S=\wedge \oplus (*,B)$ for some non-empty trees $A,B$ with $|A|=|B|=k$ for some $k \geq 1$ . By the assumption at the beginning of the proof, we have $lw(A)=\alpha$ and $rw(B)=\beta$. Without being worried about caret simplifications, a simple computation shows that a diagram for $g^n$ is $(\wedge \oplus (A_n,*),\wedge \oplus (*,B_n))$, where $A_n$ and $B_n$ can be obtained by induction as follows: $A_1=A$, $A_n=A_{n-1} \oplus^1 A$ and $B_1=B$ and $B_n=B_{n-1} \oplus^{(n-1)k+1} B$.

To finish this case, we need a technical estimation which will be also useful for the general case. We start with a definition.

\begin{definitie}
Let $T$ be any tree. \emph{The approximation from the left wing} of $T$ is the unique sequence of trees
$$T_1 \subseteq T_2 \subseteq \ldots \subseteq T_m,$$
\noindent where $T_1=T \cap L_{\infty}$ and $T_i=T_{i-1} \oplus^L$, where $L \subseteq [|T_{i-1}|+1]$ is the maximal subset of $[|T_{i-1}|+1]$ such that $T_{i-1} \oplus^L \subseteq T$.  Notice that after finitely many step (say $m$) this process of attaching carets stops and we obtain $T_m=T$. Notice that for any $i$, $T_i$ and $T_{i+1}$ belong to a same cube. $m$ is called \emph{the left approximation number of $T$}. A similar definition can be made for the right wing of $T$ (or more general by starting from any long snake).
\end{definitie}

\begin{exemplu}

If $T=L_4 \oplus$, then its approximation from the left wing is $L_5 \subseteq L_4 \oplus$ and its approximation from the right wing is $$R_2 \subseteq R_2 \oplus^1=T_2 \subseteq T_2 \oplus^{\{1,2\}}=L_3 \oplus^{\{3,4\}} \subseteq (L_3 \oplus^{\{3,4\}})\oplus^{\{1,2\}}=L_4\oplus^{\{3,4,5\}} \subseteq L_4 \oplus$$
\end{exemplu}

\begin{lema}
\label{estimare aripi}
(Distance estimation along wings) Let $T$ be a tree and let $m$  be its left approximation number. Then we have:
$$d(O,T) \leq lw(T) + \sqrt{lw(T)}(1+\sqrt{2})2^{m/2}$$
\end{lema}

\begin{remarca}
An identical estimate holds for the right wing. Notice that in general the approximation is very weak, but it will be useful in our situation when we iterate an element of $F$ and the wings of the trees start to dominate the shape of the diagram.
\end{remarca}

\noindent \emph{ Proof.} Let $$T_1 \subseteq T_2 \subseteq \ldots \subseteq T_m$$ be the approximation of $T$ starting from the left wing. By definition we have $|T_1|=lw(T)$. Notice that $|T_{i+1}| \leq 2|T_i|$ and that $d(T_i,T_{i+1})=\sqrt{|T_{i+1}-T_i|}.$ By the triangle inequality we have:
\begin{eqnarray*}
d(O,T)     &\leq& d(O,T_1)+ \sum_{i=1}^{m-1} d(T_i,T_{i+1})\\
           &\leq& lw(T) + \sum_{i=1}^{m-1}\sqrt{2^{i-1}lw(T)} \\
           &=& lw(T) + \sqrt{lw(T)}\sum_{i=1}^{m-1}\sqrt{2^{i-1}}\\
           &\leq& lw(T) + \sqrt{lw(T)}(1+\sqrt{2})2^{m/2}
\end{eqnarray*}
$\clubsuit$

\begin{figure}[!h]
\label{rolex}
\begin{center}
\includegraphics[scale=0.9]{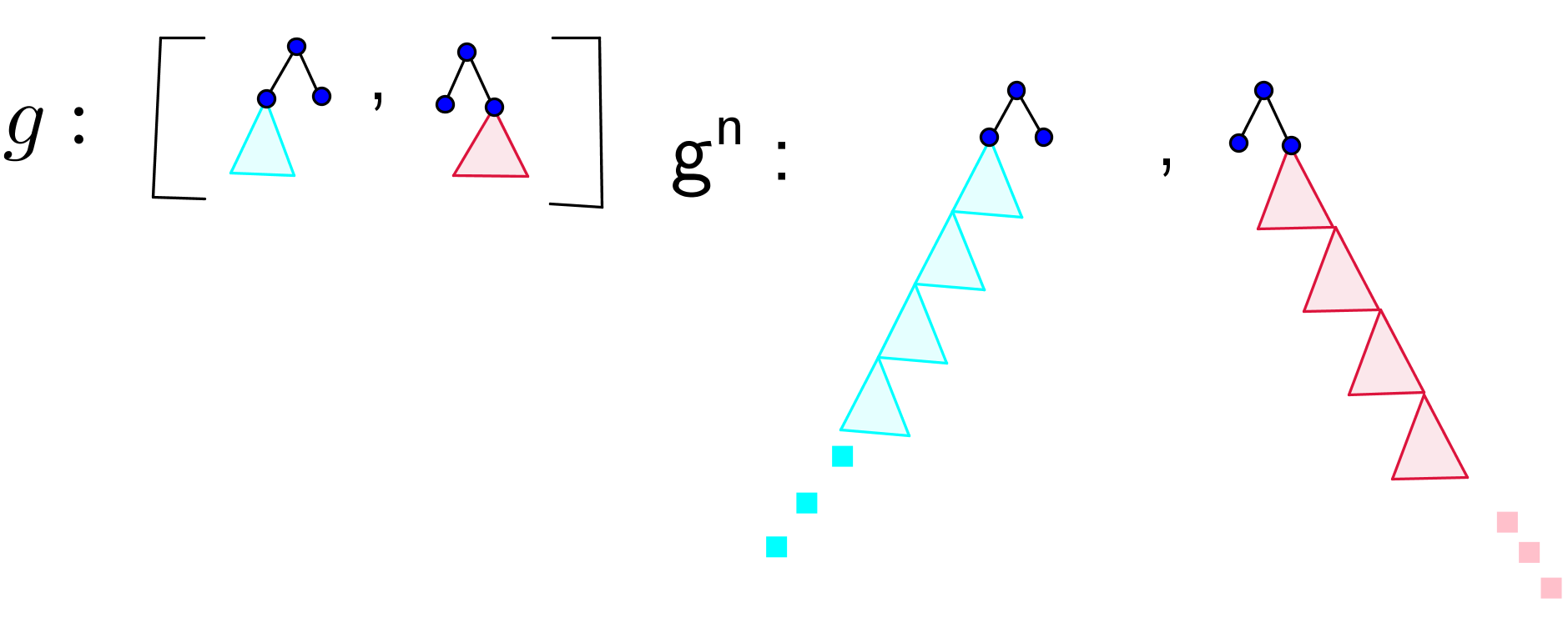}
\caption{This picture illustrates the proof of the upper bound in the case discussed so far. The two triangles colored the same in the picture represents the same binary tree. So we have a simple periodicity behavior of diagrams under iteration. The left wing of the left side grows linearly and controls the first slope. The right wing in the right binary tree also grows linearly and controls the last slope. All the other parts of the diagram grows sub-linearly: indeed, combinatorially the width of the left side doesn't grow under iteration (the same for the width of the right side); this is reflected geometrically in the Lemma 2.3.10 in the main text. For the general case of an arbitrary irreducible isometry, the situation is more complicated but the idea the same: in the periodicity six trees are involved, instead of two, and only along subsequences.}
\end{center}
\end{figure}

We can now return to our proof. By construction notice that for any $n$ the left approximation number of $A_n$ equals the left approximation number of $A$ and the right approximation number of $B_n$ equals the right approximation number of $B$. We denote these two constants $a$ and $b$ and let $K=(1+\sqrt{2})max\{2^{a/2},2^{b/2}\}$. Using \ref{calcul recursiv}, \ref{asimptotic} and \ref{estimare aripi} we have, for any $n \geq 1$:
\begin{eqnarray*}
n|g|=|g^n| &\leq& d(g^n(\wedge \oplus(A_n,*)),\wedge \oplus(A_n,*))\\
           &=& d(\wedge \oplus(A_n,*),\wedge\oplus(*,B_n)) \\
           &=& \sqrt{d(O,A_n)^2+d(O,B_n)^2}\\
           &\leq& \sqrt{[lw(A_n)+\sqrt{lw(A_n)}(1+\sqrt{2})2^{a/2}]^2+[rw(B_n)+\sqrt{rw(B_n)}(1+\sqrt{2})2^{b/2}}]^2\\
           &\leq& \sqrt{(n\alpha + K\sqrt{n\alpha})^2 + (n\beta + K\sqrt{n\beta})^2}
\end{eqnarray*}

\noindent By letting $n \rightarrow \infty$, we get $|g| \leq \sqrt{\alpha^2 + \beta^2}$.

We now return to the remaining (general) case, where if $(T,S)$ is a diagram of $g$, then $U=T \cap S \neq \wedge$, that is, $|U|=m \geq 2$. Then there are two disjoint subsets of $[m+1]$, $I$ and $J$, and some non-empty trees $A_1, \ldots, A_r$ and $B_1, \ldots, B_s$ with $r=|I| \geq 1$, $s=|J| \geq 1$, such that
$$T=U \oplus^I(A_1, \ldots, A_r), S= U \oplus^J(B_1, \ldots, B_s),$$
\noindent with $lw(A_1)=\alpha$, $rw(B_s)=\beta$, $1 \in I$ and $m+1 \in J$,  by our assumptions. We replace $g$ with another element in the same conjugacy class. Let $x=(U,L_2\oplus^2 V)$, where $V$ is any tree of size $m-2$. Without being worried about caret simplifications, a straightforward computation shows that $$x^{-1}gx=((L_2 \oplus^2 V) \oplus^I(A_1, \ldots, A_r),(L_2 \oplus^2 V)\oplus^J(B_1, \ldots, B_s))$$

We replace $g$ with $x^{-1}gx$ and we rewrite the diagram in the form $$(L_2 \oplus (X_1,X_2,*), L_2 \oplus (*,Y_2,Y_1)),$$
\noindent with $lw(X_1)=\alpha$ and $rw(Y_1)=\beta$. We can assume that $g(1/4) > 1/2$: indeed, using \ref{asimptotic}, we can replace $g$ with a sufficiently large positive power $g^m$ (notice that the diagram of $g^m$ obtained by multiplying without performing carets simplifications, has the same form as above, of course with others $X_1,X_2,Y_1,Y_2$). With the new allowed assumption we have $|X_1| > |Y_2|$ and $|X_2|< |Y_1|$. Since it is not completely trivial, we explicit the behavior of the diagrams under iteration. As a rule, we do not perform caret simplifications, we are not interested in computing the reduced diagram.

We first compute a diagram for $g^2$:
$$(L_2 \oplus (X_1,X_2,*), L_2 \oplus (*,Y_2,Y_1)) (L_2 \oplus (X_1,X_2,*), L_2 \oplus (*,Y_2,Y_1))$$
At the first stage of the computation we add carets at the right tree of the first diagram and at the left tree of the second diagram until we obtain $L_2 \oplus (*,X_2 \cup Y_2,Y_1)$ and respectively $L_2 \oplus (X_1, X_2 \cup Y_2,*)$. The corresponding carets attached  at the left tree of the first diagram are, because of our assumption, all located at the leaves of $X_1$, but not at the first leaf. We denote by $X_3$ the tree which replaced $X_1$ in this way. Similarly, for the second diagram, the corresponding carets attached at the right tree are all located at the leaves of $Y_1$, but not at the last one. We denote by $Y_3$ the tree which replaced $Y_1$ in this way. After performing this operations we get:
$$(L_2 \oplus (X_3,X_2,*),L_2 \oplus (*, X_2 \cup Y_2, Y_1)) (L_2 \oplus (X_1, X_2 \cup Y_2, *)), L_2 \oplus(*,Y_2,Y_3)) $$

We only have now to complete the middle trees with carets until they become equal with $L_2 \oplus(X_1,X_2 \cup Y_2,Y_1)$, that is we only have to attach $X_1$ at the first leaf in the first diagram and $Y_1$ at the last leaf of the second diagram. We obtain the following diagram for $g^2$:
$$( (L_2 \oplus(X_3,X_2,*)\oplus^1 X_1), (L_2 \oplus (*,Y_2,Y_3) \oplus^{-1}Y_1)),$$
\noindent where the superscript $-1$ at $\oplus$ means that we attach a tree at the last leaf. Now a diagram of $g^n$ can be immediately checked by induction to be:
$$(L_2 \oplus (A_n,X_2,*), L_2 \oplus (*,Y_2,B_n))$$
\noindent where $A_1=X_1$ and $A_n=C_{n-1} \oplus^1 X_1$, where $C_1=X_3$ and $C_k=C_{k-1} \oplus^1 X_3$. Similarly, $B_1=Y_1$ and $B_n=D_{n-1} \oplus^{-1}$ , where $D_1=Y_3$ and $D_k=D_{k-1}\oplus^{-1}$.

Notice that the left approximation number of $A_n$ is the maximum between the left approximation number of $X_1$ and $X_3$ and the right approximation number of $B_n$ is the maximum between the right approximation number of $Y_1$ and $Y_3$. Denote $a$ and $b$ these constants and let $K=(1+\sqrt{2})max\{2^{a/2},2^{b/2}\}$.

Again, using \ref{calcul recursiv}, \ref{asimptotic} and \ref{estimare aripi} we have, for any $n \geq 1$:
\begin{eqnarray*}
n|g|=|g^n| &\leq& d(g^n(L_2 \oplus(*,Y_2,B_n)),\wedge \oplus(*,Y_2,B_n))\\
           &=& d(L_2 \oplus (A_n,X_2,*),L_2\oplus(*,Y_2,B_n)) \\
           &=& \sqrt{d(O,A_n)^2 + d(X_2,Y_2)^2+ d(O,B_n)^2}\\
           &\leq& \sqrt{[lw(A_n)+\sqrt{lw(A_n)}K]^2+ d(X_2,Y_2)^2+[rw(B_n)+\sqrt{rw(A_n)}K]^2}\\
           &\leq& \sqrt{(n\alpha + K\sqrt{n\alpha})^2 + d(X_2,Y_2)^2+ (n\beta + K\sqrt{n\beta})^2}
\end{eqnarray*}

Letting $n \rightarrow \infty$, we have $|g| \leq \sqrt{\alpha^2 + \beta^2}$. The proof of \ref{lungime translatie} is now complete.

\vspace{2mm}

\begin{remarca}
It is now easy to deduce Corollary 1. Indeed, if for any $g \in F$, the translation length of $g$ is at least $\sqrt{log_2(a)^2+log_2(b)^2}$, where $a$ is the most left non-trivial slope and $b$ the most right non-trivial slope of $g$. The proof is line by line the same  with the lower bound for an irreducible element, except the fact that instead using $L_{\infty}$ and $R_{\infty}$ with the two long snakes determined by the most left and most right point of the support of $g$ as a homeomorphism of the unit interval. In fact, it is not difficult to compute all the translation lengths for the elements of $F$.
\end{remarca}

\subsection{Parabolic Isometries}

In the previous two sections we only really dealt with vertices in $X$. We give versions of the previous lemmas which applies to interior points in $X$. While the strategy is the same displayed in the section concerning the particular isometry $(L_3, R_2 \oplus^2)$, the details are more complicated. All the section is devoted to proving the second part of Theorem 1.

\begin{propozitie}
\label{parabolic}
If $g \in F$ is an irreducible element with $g'(0)g'(1) \neq 1$, then $g$ is parabolic.
\end{propozitie}

We start with the following very simple observation, but which plays a crucial role at the very end of the proof.
\begin{lema}
\label{numarare}
(Carets Counting). If $(T,T_1, \ldots, T_k)$ is a diagram representing a vertex in $X$ then we have
$$ |T|=k-1+ \sum_{i=1}^k|T_i|$$.
\end{lema}

\noindent \emph{Proof.} The number of leaves on the left side of the diagram is $|T|+1$ and the number of leaves on the right side of the diagram is $k+ \sum_{i=1}^k|T_i|$. These two numbers are equal by definition. $\clubsuit$

\begin{definitie}
\label{subdivision}
(Subdivision of a cube)

(a) Let $C=[0,1]^N$ be the euclidian cube and let $x=(t_1, \ldots, t_n)$ with $t_i \in (0,1)$ be an interior point. Two points in the cube $u=(a_1, \ldots, a_N)$ and $b=(b_1, \ldots, b_N)$ are called \emph{separated by $x$} if for any $j=1, \ldots, N$ we have $(a_i-t_i)(b_i-t_i) \leq 0$.

(b) For any $\alpha \in \{\leq,\geq \}^N$ a  list of signs $\{ \leq, \geq \}$ the subset $C_{\alpha} \subset C$ defined by $C_{\alpha}=\{ (c_1, \ldots, c_N | c_i \alpha(i) t_i$ for all $i \}$ is called \emph{a room of $C$ with respect to $x$}. A cube has $2^N$ rooms with two by two disjoint interior. Two rooms $C_{\alpha}$ and $C_{\beta}$ are called \emph{opposite} if $\alpha(i) \neq \beta(i)$ for all $i=1, \ldots, N$. Notice that two points are separated by $x$ if and only if they belong to two opposite rooms. We call two curves in $C$ separated by $x$  if their images lie entirely in opposite rooms.
\end{definitie}

\begin{lema}
\label{separation}
(Separated geodesics near an interior point) Let $M$ be a $CAT(0)$ cubical complex, $C_1,C_2$ be two cubes in $X$ with non-trivial intersection $C=C_1 \cap C_2$ and $x$ be a point in $C$. If $c_1$ and $c_2$ are two geodesics starting at $x$ with the image of $c_1$ in $C_1$ and the image of $c_2$ in $C_2$ and if the images of their projections on $C$ are separated by $x$, then $\angle_x(c_1,c_2) \geq \pi/2$.
\end{lema}

\noindent \emph{Proof.} By Lemma \ref{interplay} we can realize the sub-complex $Conv[C_1,C_2]$ in the standard integer cubulation of $\mathbb{R}^N$ for some $N$. We may assume that $C$ is embedded as $[0,1]^n$ if $n=dim(C)$ and $x$ has coordinates $(t_1, \ldots, t_n,0, \ldots, 0)$ (after identification). Fix a small $t > 0$. Let $\alpha_1, \ldots, \alpha_N$ be the coordinates of $c_1(t)$ and $(\beta_1, \ldots, \beta_N)$ the coordinates of $c_2(t)$. Since the images of $c_1(t)$ and $c_2(t)$ lie in opposite  rooms of $C$ with respect to $x$, we have $(\alpha_i-t_i)(\beta_i-t_i) \leq 0$ for all $i=1, \ldots, n$ and since $C=C_1 \cap C_2$ we also have $\alpha_j \beta_j \leq 0$ for $j=n+1, \ldots, N$. We have:

\begin{eqnarray*}
d(c_1(t),c_2(t))^2 &\geq& d_N(c_1(t),c_2(t)^2 \\
                   &=& \sum_{k=1}^N(\alpha_k-\beta_k)^2 \\
                  &=& \sum_{i=1}^n[(\alpha_i-t_i)-(\beta_i-t_i)]^2+ \sum_{j=n+1}^N(\alpha_j-\beta_j)^2\\
                  &\geq& \sum_{i=1}^n(\alpha_i-t_i)^2+\sum_{j=n+1}^N\alpha_j^2 + \sum_{i=1}^n(\beta_i-t_i)^2+\sum_{j=n+1}^N\beta_j^2\\
                  &=& 2t^2
                  \end{eqnarray*}

\noindent where $d_N$ is the euclidian distance on $\mathbb{R}^n$. Letting $t \rightarrow 0$ we are done. $\clubsuit$

\vspace{2mm}

\begin{definitie}
\label{proiectie vizuala}
Given a point $x \in X$ with reduced diagram $$(T,T_1, \ldots, T_k)[t_1, \ldots, t_k]$$  ( $0 \leq t_i < 1$), the point $T[s_1, \ldots, s_r] \in \mathcal{T}$ is called the \emph{visual projection} of $x$, where the $s_i$'s are defined as follows ($r=|T|+1$). Let $I \subset [k]$ be the set of indices for which $T_i=*$ and let $\sigma: [r] \rightarrow [k]$ defined by $\sigma(i)=j$ if the $i$th leaf of the right-side part of the diagram  $(T,T_1, \ldots,T_k)$ belongs to $T_j$. We define $s_i=t_{\sigma(i)}$ if $\sigma(i) \in I$ and 0 otherwise. We denote this point by $\Pi(x)$. For example, if $x=(L_3,*, \wedge,*)[1/3,1/4,0]$, then $\Pi(x)=L_3[1/3,0,0,0]$.
\end{definitie}

\begin{lema}
\label{proiectie copaci}
(Projection on Trees) If $x \in X$ is a point in the complex, then the projection of $x$ on the sub-complex $\mathcal{T}$ is $\Pi(x)$, where $\Pi(x)$ is the visual projection of $x$.
\end{lema}

\noindent \emph{Proof.} We keep the notations from the previous definition, $T=(T_1, \ldots, T_k)$, $\Pi(x)=(T,s_1, \ldots, s_r)$. We show that $\angle_T(c,c') \geq \pi/2$, where $c_1$ is the geodesic from $\Pi(x)$ to $x$ and $c_2$ is any geodesic starting at $\Pi(x)$ with the image in $\mathcal{T}$. Fix a sufficiently  small $t >0$.

Notice that $x$ and $\Pi(x)$ belong to $Conv[T \oplus^J,(T,T_1, \ldots, T_k)]$, where $J=\{j \in [r] | s_j > 0\}$. As a consequence, by \ref{interplay}, $c_1(t)$ belongs to the cube $C_1=C(T, \Delta_1| K)$, where $\Delta_1$ is the termination of the sequence $(T_1, \ldots, T_k)$ (see \ref{terminatie}) and $K$ consist of the positions where we have a $\wedge$ and the positions where we have $*$ with the leaf labeled by a number in $J$. By \ref{interplay} again, $c_2(t)$ belongs to a cube $C_2=C(T, \Delta_2 |L)$, where $\Delta_2$ is a sequence of trees containing only $*$ and $\wedge$, the $\wedge$ being labeled with pairs of leaves which correspond to carets in $T$ (otherwise we are not in $Trees$ anymore) and $L$ contains all the positions where we have $\wedge$ and other positions where we have $*$ included all those labeled by a number from $J$. Notice also that $(T, \Delta_2)$ is not necessarily reduced. Nevertheless, since the diagram $(T,T_1, \ldots, T_n)$ is reduced, we have that $C=C_1 \cap C_2=C=C(T|J)$.

We want to apply \ref{separation} to the configuration $C=C_1 \cap C_2$. To prove the separation of $c_1(t)$ and $c_2(t)$ (see \ref{subdivision}), it is enough to show that the coordinate of $c_1(t)$ at any position $j \in J$ is $s_j=t_{\sigma(j)}$ ($\sigma$ is the map defined in \ref{proiectie vizuala})(since in this case we can choose $c_1(t)$ to belong in any room, in particular in the opposite room where $c_2(t)$ belongs). To do this let us notice that the sub-complex $Z=Conv[T \oplus^J,(T,T_1, \ldots, T_k)]$ splits as $Y \times [0,1]^J$, where $Y=Conv[T,(T_1, \ldots, T_k)]$. Indeed, to see this we write all the vertices in $Z$ in the reduced form $(T \oplus^{J_1},\Delta)$, where $\Delta$ denotes (by abuse of language) only the configuration of trees which involves the positions labeled initially with numbers outside $J$, the operations performed at the positions with leaves labeled initially in $J$  being recorded on the left by the subset $J_1 \subseteq J$. Since it is easily seen that the operations at positions with the leaf labeled in $J$ required to reach  $T \oplus^J$ by moving from one vertex to another from $(T,T_1, \ldots, T_k)$ are independent from those involving leaves labeled outside $J$ we have that the map $\phi: Z \rightarrow Y \times [0,1]^J$ defined by :

$$   \phi(T \oplus^{J_1}, \Delta)=((T, \Delta), \chi(J_1))$$

\noindent where $\chi$ is the characteristic function of $J$, is  cubical and bijective (and hence an isometry).

Since both $\Pi(x)$ and $x$ have the same coordinates at $J$ (namely $s_j=t_{\sigma(j)}$ for any $j \in J$), it now follows that the geodesic from $\Pi(x)$ to $x$ has constant coordinates at these positions and we are done. $\clubsuit$

\vspace{2mm}

\begin{lema}
\label{proiectie serpi int}
(Projection on snakes) Let $x=T[t_1, \ldots, t_k]$ be a point in $\tau$ and let $\sigma$ be a long snake. The projection of $x$ on the half-line geodesic associated with $\sigma$ is $(T \cap \sigma)[0, \ldots,t_j, \ldots,0]$, where $j$ is the label of the leaf in $T$ at the vertex connecting $T \cap \sigma$ with the next caret in $\sigma$.
\end{lema}

\vspace{2mm}

\noindent \emph{Proof.} The proof is very similar with \ref{proiectie serpi}, here in addition, we have to locate the rooms in which the geodesics starts in order to apply \ref{separation}. Let $n$ be such that $T \cap \sigma=\sigma[n]$. We denote $z=\sigma[n][0, \ldots,t_j, \ldots,0]$, $c$ the geodesic from $z$ to $x$, $c_1$ the geodesic from $z$ to $\sigma[n]$ and $c_2$ the geodesic from $z$ to $\sigma[n+1]$. We need to show that $\angle_z(c,c_1) \geq \pi/2$ and $\angle_z(c,c_2) \geq \pi/2$. We focus on the first angle. We can assume $0 < t_j < 1$, otherwise we can apply \ref{proiectie serpi} directly.

By \ref{interplay}, the geodesic $c$ starts in the cube $C=C(\sigma[n] \oplus I)$, where $I$ is a subset of $[n+1]$ containing $k$, the label of the leaf at which $\sigma[n]$ connects to $\sigma[n+1]$. The geodesic $c_1$ starts, of course, in the cube $C_1=C(\sigma[n]|\{k\})$. We have $C_1 \cap C=C_1$. For $c$ we have both choices of room in $C$ and for $c_1$ we have only the choice $\leq$. So by choosing the sign $\geq$ for $c$ we are done using \ref{separation}.

The second inequality follows the same way, with the only difference that we have to choose the room defined by $\geq$ for $c_2$. $\clubsuit$

\vspace{2mm}

\noindent \emph{Proof of Proposition \ref{parabolic}.} Let $g \in F$ be an irreducible element with $g'(0)g'(1) \neq 1$. Replacing $g$ with $g^{-1}$, if necessary, we can assume that $\alpha=log(g'(0)) > 0$ and $\beta=-log(g'(1)) > 0$. Let $x$ be any point in $X$ and consider $(A,A_1, \ldots, A_k)[t_1, \ldots, t_k]$ to be a diagram representing $x$ with all $t_i < 1$. The goal is to prove that $d(gx,x) > \sqrt{\alpha^2 + \beta^2}$. For the moment we assume $k \geq 2$. The case $k=1$ is easy but somewhat doesn't fit in our general proof, and so we postpone it until the end.

About the trees $A$ and $A_i$ we may assume that we added enough carets at the reduced diagrams of $x$ and $g$ such that a diagram of $g$ is of the form $(B,A)$, for some tree $B$. Also denote by $\overline{x}$ the vertex of $X$ represented by the diagram $(A,A_1, \ldots, A_k)$ and by $\gamma$ and $\gamma'$ the binary logarithms of the first and the last slope of $\overline{x}$ (when viewed as a homeomorphism). We have then $lw(B)=lw(A)+ \alpha$, $ rw(B)=rw(A)+ \beta$, $lw(A)=lw(A_1)+ \gamma$ and $rw(A)=rw(A_k)+\gamma'$.

\vspace{2mm}
\noindent \emph{Step1. Reducing the computation.}
\vspace{2mm}

We will in fact show that $d(hgx,hx) > \sqrt{\alpha^2 + \beta ^2}$, where $h \in F$ is some convenient isometry. Indeed, let $h \in F$ be the element represented by the diagram $(L_{k-1} \oplus(A_1, \ldots,A_k),A)$. Notice that $log(h'(0))=k-1-\gamma$ and $log(h'(1))=1- \gamma'$. We have that $hx$ is represented by the diagram $L_{k-1}[t_1, \ldots, t_k]$. We cannot say very much about the diagram of $hgx$, but it will be enough to notice that $log((hg\overline{x})'(0))=\alpha+k-1$ and $log((hg\overline{x})'(1))=1-\beta$.

The proof is now reduced to show that $d(M,N) > \sqrt{\alpha^2 + \beta^2}$, where $N$ is of the form $L_{k-1}[t_1, \ldots, t_k]$ (with $0 \leq t_i < 1$) and the reduced diagram of $M$ is of the form $(T,T_1, \ldots, T_k)[t_1, \ldots,t_k]$ with $lw(T)=lw(T_1)+ \alpha+k-1 $ and $rw(T_k)=rw(T)+ \beta-1$. So from now on we focus only on points of this type.

\vspace{2mm}
\noindent\emph{Step2. Breaking the diagram, the obtuse triangle}
\vspace{2mm}

Let $M=(T,T_1, \ldots, T_k)[t_1, \ldots,t_k]$ and $N=L_{k-1}[t_1, \ldots, t_k]$ be like at the end of the previous step. Let $P=\Pi(M)=T[s_1, \ldots, s_r]$ the projection of $M$ on the sub-complex $\mathcal{T}$, like in \ref{proiectie vizuala} and \ref{proiectie copaci}. By the same Lemma, the triangle $M,P,N$ is obtuse at $P$, so we have $d(M,N)^2 \geq d(M,P)^2 + d(P,N)^2$. If $h' \in F$ is the element represented by the diagram $(L_{k-1} \oplus (T_1, \ldots, T_k), T)$ (notice that by \ref{numarare}, the element is well defined), then we have $d(P,M)=d(h'P,h'M)$, which in the diagram language can be written

$$d((L_{k-1} \oplus (T_1, \ldots, T_k))[s_1, \ldots, s_r], L_{k-1}[t_1, \ldots, t_k]).$$

Using \ref{calcul recursiv} and projecting on $L_{\infty}$ and $R_{\infty}$, we have:
\begin{eqnarray*}
d(M,N)^2 &\geq& d(N,P)^2 + d(h'M,h'P)^2 \\
                   &\geq& d(N,P)^2 + d(\phi_{L_{k-1}}^{-1}(h'M),\phi_{L_{k-1}}^{-1}(h'P))^2 \\
                  &\geq& [(lw(T)+s_1)-(k-1+t_1)]^2 + [(rw(T_k)+s_r)-t_k]^2
                  \end{eqnarray*}

\vspace{2mm}

\noindent \emph{Step3. Ruling out most of the configurations}

\vspace{2mm}

During this step of the proof, we assume $rw(T) \geq 2$, that is $rw(T_k) \geq \beta+1$.

\vspace{2mm}

Case1. $lw(T) \geq k + \alpha$. This implies $lw(T_1)>0$ and thus $s_1=0$. The inequality in the previous step becomes:

$$ d(M,N)^2 \geq (\alpha+1-t_1)^2 + (\beta + 1 - t_k)^2 > \alpha^2 + \beta^2$$

Case2. $lw(T)=k-1+ \alpha$. This implies $T_1=*$ and $s_1=t_1$. The inequality in the previous step becomes:

$$d(M,N)^2 \geq \alpha^2 + (\beta+1-t_k)^2 > \alpha^2 + \beta^2$$

From now on we can assume that $rw(T)=1$, that is $rw(T_k)=\beta \geq 1$.

\vspace{2mm}

\noindent \emph{Step4. Reconfiguration, perturbing the obtuse triangle}

\vspace{2mm}

 The only case left is when  $M=(T, T_1, \ldots, T_k)[t_1, \ldots, t_k]$ and $N=L_{k-1}[t_1, \ldots, t_k]$ with $rw(T)=1$, that is $rw(T_k)=\beta$. In this situation we break the diagram in a slightly different way. Recall that $P=\Pi(M)=T[s_1, \ldots, s_r]$ like in \ref{proiectie vizuala}. In our case $s_r=0$ since $T_k$ is non-empty. We consider the point $P_1=T[s_1, \ldots, s_{r-1},t_k]$ which is exactly like $P$ with the exception of the last coordinate, which is now equal with $t_k$.  We claim the following fact, which will be proved in the next step.

\vspace{2mm}
\noindent \emph{Claim.} The triangle $M,P_1,N$ is obtuse at $P_1$.
\vspace{2mm}

We have then:
$$d(M,N)^2 \geq d(N,P_1)^2 + d(M,P_1)^2 = d(N,P_1)^2 + d(h'M,h'P_1)^2$$
\noindent where $h'=(L_{k-1} \oplus (T_1, \ldots, T_k), T) \in F$ is like in the previous step.

If $lw(T) \geq \alpha+k$ then $T_1 \neq *$ and $s_1=0$. We have then

$$d(N,P_1) \geq d(\Pi_{L_{\infty}}(P_1), \Pi_{L_{\infty}}(N))=lw(T)+s_1-(k-1+t_1) \geq \alpha+1-t_1$$
and
$$d(h'M,h'P_1) \geq d(\Pi_{R_{\infty}}{h'M}, \Pi_{R_{\infty}}(h'P_1))=(1+rw(T_k)+t_k)-(1+t_k)=rw(T_k)= \beta$$
In this case we have  the desired inequality and thus we can also assume from now on that $lw(T)=\alpha+k-1$, that is $T_1=*$ and $s_1=t_1$. Recall also that we are in a stage of a proof where we assume $rw(T_k)=\beta$, that is $rw(T)=1$.

We now apply the counting in $\ref{numarare}$ to the diagram $(T,T_1, \ldots, T_k)$:
$$|T|=k-1+|T_1|+ \ldots +|T_k|$$

\vspace{2mm}

Case1. If $\alpha > \beta$, then from the previous counting and assumption there must be a caret belonging to one of trees $T_2, \ldots, T_{k-1}$ (recall that $T_1$ is empty by now) or to $T_k$ but in this case not to its right wing. If the caret belongs to some $T_j$ with $2 \leq j \leq k-1$ then $T_j$ is non-empty and  by \ref{calcul recursiv} we have:

\begin{eqnarray*}
d(M,N)^2 &\geq& d(N,P_1)^2 + d(h'M,h'P_1)^2 \\
                &\geq& \alpha^2 + d^2(*[t_j],T_j) + \beta^2 \\
                &\geq& \alpha^2 + \beta^2 + (1-t_j)^2
\end{eqnarray*}

Since all the $t_i$'s were chosen to be $ < 1$ we are done in this case.

If the caret is located in $T_k$, then we have

\begin{eqnarray*}
d(M,N)^2 &\geq& d(N,P_1)^2 + d(h'M,h'P_1)^2 \\
         &\geq& \alpha^2 + [d(*[t_k],T_k[0,\ldots,0,t_k])]^2 \\
         &\geq& \alpha^2 + d(*[t_k],R_{\beta}[0, \ldots, t_k])^2 + d(R_{\beta}[0, \ldots, t_k],T_k[0, \ldots, 0,t_k])^2\\
         &>& \alpha^2 + \beta^2.
\end{eqnarray*}

\noindent The last two inequalities follows from the fact that the triangle
$$T_k[0, \ldots,0, t_k], R_{\beta}[0, \ldots,0,t_k],*[t_k]$$
is obtuse at the secondly named vertex, which is just the projection of the first one on $R_{\infty}$. Because of the presence of the extra caret, $T_k[0, \ldots,0,t_k]$ is not situated on $R_{\infty}$ and hence the triangle is not degenerated.

\vspace{2mm}

Case2. If $\alpha < \beta$, then the previous counting shows that there must be a caret in $T$ not situated on the left wing of $T$. So we have either the presence of a caret in $T$ glued at the $j^{th}$ leaf of $L_{k-1}$ for some $j \geq 2$ or a caret in the subtree of $T$ rooted at the first leaf of $L_{k-1}$ but not on its left wing. If we are in the first situation then we have:
\begin{eqnarray*}
d(M,N)^2 &\geq& d(N,P_1)^2 + d(h'M,h'P_1)^2 \\
         &\geq& \alpha^2 + (1-t_j)^2 + \beta^2
\end{eqnarray*}

\noindent and if we are in the second situation and we denote by $T'$ the subtree of $T$ rooted at the first leaf of $L_{k-1}$, then, keeping in mind that $lw(T')=\alpha$, we have
\begin{eqnarray*}
d(M,N)^2 &\geq& d(N,P_1)^2 + d(h'M,h'P_1)^2 \\
                &\geq& d(*[t_1],T'[t_1, \ldots])^2 + \beta^2 \\
                &\geq& d(*[t_1], L_{\alpha}[t_1, 0, \ldots, 0])^2 + d(L_{\alpha}[t_1, 0, \ldots, 0],T'[t_1,  \ldots])^2 + \beta^2 \\
               &>& \alpha^2 + \beta^2
\end{eqnarray*}

\noindent The last two inequalities follows from the fact that the triangle
$$T'[t_1, \ldots], L_{\alpha}[t_1, \ldots,0,0],*[t_k]$$
is obtuse at the secondly named vertex, which is just the projection of the first one on $L_{\infty}$. Because of the presence of the extra caret, $T'[t_1, \ldots]$ is not situated on $L_{\infty}$ and hence the triangle is not degenerated.

At this point we only have to check the claim made about the triangle $M,P_1,N$ and to come back and solve the case $k=1$.

\vspace{2mm}

\noindent \emph{ Step 5. Proof of the claim made at Step 4}

\vspace{2mm}

We show that the triangle $M,P_1,N$ is obtuse at $P_1$. The proof is very similar with \ref{proiectie copaci} to which we refer the reader. Here we just record the variations needed to conclude. Let $t>0$ be sufficiently small, let $c_1$ be the geodesic from $P_1$ to $M$ and $c_2$ be the geodesic from $P_1$ to $N$. The two geodesics start in two cubes $C_1$ and $C_2$ defined like in \ref{proiectie copaci} and which intersect in $C=C_1 \cap C_2=C(T|J)$. Nevertheless there are differences in the rooms they occupy. We stress these differences.

For the geodesic $c_1$ notice if $t_k=0$ there is no difference: $P_1=P=\Pi(M)$, but if $t_k >0$ then $r \in J$ and because of this it is more convenient to split $Z=Conv[T \oplus^J,(T,T_1, \ldots,T_k)]$ first as $Y_1 \times Y_2$, where $Y_1=Conv[T \oplus^{ J-\{r\}},(T,T_1, \ldots, T_{k-1},*, \ldots, *)]$ and $Y_2=Conv[T \oplus^r,\\ (T,*, \ldots*,T_k)]$. The splitting holds since all the operations performed on $(T_1, \ldots, T_{k-1})$ are independent of those performed on $T_k$, when moving on the median interval from  $(T,T_1, \ldots, T_k)$  to $T \oplus^J$. $Y_1$ further splits in $Conv[T,(T,T_1, \ldots, T_{k-1}, *, \ldots,*)] \times [0,1]^{J-\{r\}}$ like in \ref{proiectie copaci}. This splitting secures all the coordinates of $c_1(t)$ at the positions labeled by $J$ to be constant ($=s_j=t_{\sigma(j)}$), with the exception of the coordinate at the leaf labeled by $r$ (in fact it is clear that the sign choice for this coordinate will be $\geq$, but we don't need a precise choice). So the image of $c_1(t)$ in the cube $C=C(T|J)$  belong to all the rooms (with respect to $P_1$) with a forced sign choice just for the last coordinate.

We now determine to which rooms $c_2(t)$ belong. To do this recall that $rw(T)=1$. The image of $c_2$ is in $Z= Conv[L_{k-1},T \oplus^J]$, which by a same splitting argument is $Y \times [0,1]$, where $Y=Conv[L_{k-1},T \oplus^{J-\{r\}}]$. Since $M$ and $P_1$ have the same last coordinate ($=t_k$), it follows that the last coordinate is constant along $c_2(t)$ also, so we can choose at least two rooms in which $c_2(t)$ belongs, namely with $\geq$ and $\leq$ at the last coordinate.

Combining the last two paragraph we can choose $c_1(t)$ and $c_2(t)$ in opposite rooms. Indeed, for the first $|J|-1$ coordinates we can choose any sign for $c_1$ and for the last one we can choose any sign for $c_2$. By \ref{separation} we have $\angle_{P_1}(c_1,c_2) \geq \pi/2$.

\vspace{2mm}

\noindent \emph{Step6. The one dimensional case.}

\vspace{3mm}

We return to the case $k=1$ which we ignored so far to avoid confusion (since we worked at both first and last coordinate which coincide in this case). In this case $x=(A,B)[t]$ where $\overline{x}=(A,B)$ is the diagram of an element in the Thompson group. The reduction in Step2 is just the conjugation of $g$ with $\overline{x}$ in this case. So the problem is reduced to show $d((T,S)[t],*[t])^2 > \sqrt{\alpha^2+\beta^2}$, where $(T,S)$ is a reduced diagram of an element in $F$ with the logarithm of the first slope $\alpha$  and the logarithm of the last slope $-\beta$, i.e. $lw(T)=lw(S)+\alpha$ and $rw(S)=rw(T)+\beta$. Since the triangle $(T,S)[t],T,*[t]$ is obtuse, by projecting on $L_{\infty}$ and $R_{\infty}$, we have:

\begin{eqnarray*}
d((T,S)[t],*[t])^2 &\geq& d(T,*[t])^2 + d((T,S)[t],T)^2 \\
                &\geq& d(T,*[t])^2 + d((S,T)(T,S)[t],(S,T)T)^2 \\
                &=& d(T,*[t])^2 + d(*[t],S)^2\\
                &\geq& (lw(T)-t)^2 + d(rw(S)-t)^2\\
                &\geq& (\alpha+1-t)^2 + (\beta+1-t)^2 \\
               &>& \alpha^2 + \beta^2.
\end{eqnarray*}

$\clubsuit$

\newpage

\section{On the Tits Boundary}

In this chapter, we discuss some aspects concerning the Tits Boundary. We only focus on the sub-complex $\mathcal{T}$. Recall from the section 1.4. that the vertices of this complex are all the elements in $Trees$ and that two trees $T,S$ are neighbors in the 1-skeleton if and only if $T$ can be obtained from $S$ by either adding or removing a caret. So a tree $T$ has $|T|+||T||+1$ neighbors, where recall that $||T||$ is the number of free carets of $T$ (and $|T|$, as always, is the number of total carets). In section 3.1. we explain Corollary 2 from the Introduction, and in section 3.3, we prove Theorem 2 and Corollary 3. Section 3.1 and 3.2 are only to orient the reader, they follow from the result in section 3.3.

\subsection{Snakes in Trees}

We first describe the median structure of $\mathcal{T}$. By  Chepoi's work in \cite{chepoi} it is known that the skeleton of any $CAT(0)$ cubical complex is a median graph, and conversely any median graph can be completed to a $CAT(0)$ cubical structure if we fill in all the maximal (Hamming) cubes. We denote by $\rho$ the median distance on $Trees$. We will use the term \emph{move} when moving from one vertex to an adjacent one.

\begin{lema}
For any $T,S \in Trees$ we have $\rho(T,S)=|T \Delta S|$.
\end{lema}

\noindent \emph{Proof.} First, we notice that $\rho(T,S) \geq |T \Delta S|$, since for any path from $T$ to $S$, at each move we either remove or add a caret.
We show that we can always find a path of length $|T \Delta S|$. We construct a path from $T$ to $S \cup T$ and one from $S$ to $S \cup T$ and then we take their concatenation. Put the obvious lexicographic ordering on $Bin$ (see \ref{def tree} and run over this set in order. Every time when we find an element of $S-T$ we add a caret and so we add a new vertex in the path (by induction we can do this while keeping the new set in Trees). This requires $|S-T|$ steps. Doing the same for $S$ we are done. $\clubsuit$

\begin{lema}
(The median structure on $Trees$) If $T,S,U \in Trees$, then $M=(T \cap S) \cup (T \cap U) \cup (S \cap U)$ is the median point of $T,S$ and $U$.
\end{lema}

\noindent \emph{Proof.} It is easy to check that $\rho(T,U)+ \rho(S,U)= \rho(S,T)+ 2(|(T\cap S)-U|+|U-(T \cup S)|)$. So the interval with endpoints $T$ and $S$ with respect to this metric is $I(T,S)=\{ U \in Trees | T \cap S \subseteq U \subseteq S \cup T\}$. If $W \in I(T,S) \cap I(T,U) \cap I(S,U)$, then we must have $(T \cap S) \cup (T \cap U) \cup (S \cap U) \subseteq W \subseteq T \cup S, T \cup U, S \cup U$. The only element of $Trees$ which satisfies this relation is $M$. $\clubsuit$

\vspace{2mm}

From now on we work our way to the Corollary 2 from the Introduction.

\vspace{2mm}

We first classify the combinatorial hyperplanes of $\mathcal{T}$. Recall that a \emph{snake} is a tree with exactly one free caret (and hence without blocked carets). The set of all snakes is denoted by \emph{Snakes}.
We specify an edge in Trees by a pair $(T,(k,l))$ where $T$ is a tree and $(k,l) \in Bin$ is a free caret of $T$; this notation represents the edge $(T,T')$, where $T'$ is the tree obtained from $T$ by deleting the caret $(k,l)$. Notice that the snakes are one to one with $Bin$ (the positions of carets), since for each $(k,l) \in Bin$ we can construct a unique snake with the free caret $(k,l)$. Also for snakes it is unnecessary to mark the unique free caret. In the next proposition equivalence means that they belong to the same (combinatorial) hyperplane.

\begin{lema}
(Hyperplanes) Any edge in $Trees$ is equivalent with an edge indexed by a snake. Every two distinct edges in $Snakes$ are not equivalent. Moreover, for any $S \in Snakes$ we have that $[S]$ (the hyperplane corresponding to the edge indexed by $S$) consists of all edges $(T,(k,l))$, where $(k,l) \in Bin$ is the free caret of $S$.
\end{lema}

\emph{Proof.} We prove the first claim. Let $(T,(k,l))$ with $T \in Trees$ and $(k,l) \in Bin$ be an edge. If $T$ is a snake we are done, if not, pick another free caret of $T$ say $(k',l') \in Bin$. Let $T'$ be the tree obtained from $T$ by deleting the caret $(k',l')$. Then the edges $(T',(k,l))$ and $(T,(k,l))$ are opposite in the square with opposite vertices $T$ and $U$, where $U$ is the tree obtained by deleting both $(k,l)$ and $(k',l')$. Therefore $(T',(k,l))$ and $(T,(k,l))$ are equivalent. By continuing in the same way, we obtain in the end that $(T,(k,l))$ is equivalent with a snake for which $(k,l)$ is the free caret.

For the second and third assertion, let us notice that if $(T,(k,l))$ is an edge, then all the edges which are square opposites to it are indexed by a tree for which $(k,l)$ is free and marked. By the definition of a hyperplane the conclusion follows. $\clubsuit$

\vspace{2mm}

We fix from now on the origin of $\mathcal{T}$ to be the empty tree $*$ and for any hyperplane $H$ we set $H^+$ to be the half-space determined by $H$ not containing the origin. Sometimes we will also denote by $H^{+}$ the set of vertices in $\mathcal{T}$ which lie in the positive half-space determined by $H$. If $S$ is a snake it is clear that the vertices in $S^{+}$ are those $T \in Trees$ with the property that $S \subseteq T$. It is also clear that the geometric hyperplane corresponding to a snake $S$ consist of all the points $T[t_1, \ldots, t_n]$ (written in reduced form) such that the head of the unique free caret of $S$ is a leaf when viewed in $T$ and the coordinate at that leaf is $1/2$.

We now study the profiles of $\mathcal{T}$. The profiles are a combinatorial approximation of the boundary of a $CAT(0)$ cubical complex, introduced by Farley \cite{farley08}, section 2.4. A profile is a collection of positive half-spaces satisfying some axioms such that it is likely that a geodesic will travel exactly in these half-spaces.

\begin{definitie}
(Farley) Let $M$ be a locally finite $CAT(0)$ cubical complex with a fixed origin $O$.

(a) Let $\leq$ be the following partial order on the hyperplanes of $M$: if $H_1,H_2$ are two hyperplanes, then $H_1 \leq H_2$ if $H_2^+ \subseteq H_1^+$. We  can also view this poset on the set of positive half-spaces, namely, $H_1^+ \leq H_2^+$ if and only if $H_1 \leq H_2$.

(b) A \emph{profile} $\mathcal{P}$ is a collection of hyperplanes satisfying the following three properties:

(i) For any finite subset $\{H_1, \ldots, H_n\}$ of $\mathcal{P}$, we have $H_1^+ \cap \ldots \cap H_n^+ \neq \emptyset$.

(ii) The poset $(\mathcal{P},\leq)$ has no maximal elements.

(iii) If $H_1 \in \mathcal{P}$ and $H_2 \leq H_1$, then $H_2 \in \mathcal{P}$.
\end{definitie}

\begin{exemplu}
Given a geodesic ray starting at the origin, the set of all half-spaces crossed by the image of the geodesic form a profile (\cite{farley08}).
\end{exemplu}

\begin{remarca}
In our case the poset relation is very simple to understand: if $S_1$ and $S_2$ are two snakes, then $S_1^+ \leq S_2^+$ if and only if $S_2^+ \subseteq S_1^+$ if and only if $S_1 \subseteq S_2$.
\end{remarca}

The profiles in $\mathcal{T}$ are classified by the trees we call closed.

\begin{definitie}
An infinite tree $\tau$ is called \emph{closed} if none of its carets is free.
\end{definitie}

For a closed tree $\tau$ we denote by $\mathcal{P}(\tau)$ the set of all snakes (hyperplanes) $S$ with the property that there is a long snake $\sigma \subseteq \tau$ and a positive integer $n$ such that $S=\sigma[n]$. In other words, we take all the long snakes included in $\tau$ and then all the possible truncations.

\begin{lema}
For any closed tree $\tau$, $\mathcal{P}(\tau)$ is a profile in $\mathcal{T}$ and any profile of $\mathcal{T}$ is of this form.
\end{lema}

\noindent \emph{Proof.} Let $\tau$ be a closed tree. We show that $\mathcal{P}(\tau)$ is a profile. We check the axioms:

(i) Given $n$ snakes $S_1, \ldots, S_n$, $S_1^+ \cap \ldots \cap S_n^+$ contains the tree $S_1 \cup \ldots \cup S_n$, and hence is non-empty.

(ii) The poset $(\mathcal{P})(\tau),\leq)$ has no maximal elements. Indeed, if $S$ is a snake in $\mathcal{P})(\tau)$ then, since $\tau$ is closed, there is a caret in $\tau$ attached either at the left or right son if the unique caret of $S$. If we call $S'$ the snake obtained by attaching this caret at $S$ we have $S \subseteq S' \subseteq \tau$.

(iii) If we have a snake $S_2 \in \mathcal{P}(\tau)$ and $S_1 \leq S_2$ is another snake, and if $S_2=\sigma[n]$ for some long snake $\sigma \subseteq \tau$ and $n$ a positive integer, then $S_1=\sigma[m]$ for some $m \leq n$ and we are done.

We now consider an arbitrary profile $\mathcal{P}$ (that is, a collection of snakes) in $\mathcal{T}$. For any $S \in \mathcal{P}$, we show there is a long snake $\sigma$ which includes $S$ and each of its truncation is in $\mathcal{P}$. Indeed, using axiom (ii), there is a snake $S' \in \mathcal{P}$ containing $S$. By (iii) all the snakes included in $S'$  also belong to $\mathcal{P}$. Replacing $S$ with $S'$ and continuing by induction we obtain a snake with the desired properties. Consider now all the long snakes obtained by this procedure starting with any snake $S \in \mathcal{P}$. The tree $\tau$ defined as the union of all these long snakes is closed (like any union of long snakes). We have $\mathcal{P}=\mathcal{P}(\tau)$.  $\clubsuit$

\vspace{3mm}

All the results in the rest of this section are particular cases of the section 3.3. The reader can skip them in principle. For us the $n$-dimensional positive orthant of the euclidian space $\mathbb{R}^n$ is the space $\{(x_1, \ldots, x_n) | x_1, \ldots, x_n \geq 0 \} \subseteq \mathbb{R}^n$. The plan is to show that the finite families of snakes generate isometrically embedded positive orthants in $\mathcal{T}$ and then to upgrade them at infinity.

Recall from \ref{serpi} that each long snake $\tau$ induces a geodesic half-line denoted $c_{\tau}$. We denote $\xi(\tau)$ the points induced on the $CAT(0)$ boundary by the geodesics $c_{\tau}$. We call them \emph{special points}. The profiles of the special points are the one corresponding long snakes (which are closed).

We move on to construct the orthants. Let $T$ be a tree, let $I$ be a subset of $[|T|+1]$ with $n$ elements $i_1 < \ldots < i_n$ and let $\tau_1, \ldots , \tau_n$ be a family of long snakes. Consider the following subcomplex of $\mathcal{T}$, denoted by $F(T, I, \tau_1, \ldots, \tau_n)$ : the vertices are $\{T \oplus^I (\tau_1[k_1], \ldots \tau_n[k_n] | k_1, \ldots, k_n \geq 0 \}$ and the cubes are all the cubes of $\mathcal{T}$ with vertices in this set. Notice that all the maximal cubes are of dimension $n$ having (in our notation) the vertex at $(0, \ldots, 0)$ coordinate $T \oplus^I (\tau_1[k_1], \ldots ,\tau_n[k_n])$ and the one at $(1, \ldots, 1)$ coordinate $T \oplus^I (\tau_1[k_1+1], \ldots ,\tau_n[k_n+1])$.

\begin{lema}
$F(T, I, \tau_1, \ldots, \tau_n)$ is an orthant of dimension $n$.
\end{lema}

\noindent \emph{Proof.} Consider the lantern map $\phi_T$, defined just before \ref{calcul recursiv} . The set $F(T, I, \tau_1, \ldots, \tau_n)$ is the image under this map of a direct product of $n$ half-lines geodesics (\ref{serpi}), so the conclusion follows. $\clubsuit$

\vspace{2mm}

An orthant of this type will be called \emph{special orthant}.

For a long snake $\tau$ and a positive integer $k$ we denote by $\tau[[k]]$ the snake obtained by removing the first $k$ carets (from the origin).

Denote by $\Game$ the set of all long snakes. Let $\tau \in \Game$ and let $T \in Trees$. We are interested in finding the half-line starting at $T$ and asymptotic with $c_{\tau}$. Notice that $T \cap \tau=\tau[r]$ for some $r \geq 0$. Define $C_{\tau,T}:[0, \infty] \rightarrow \mathcal{T}$ by $C_{\tau,T}(n)=T \cup \tau[n+r]$ and completed by  geodesic pieces. $C_{\tau,T}$ is the image of a copy of the half-line geodesic $C_{\tau[[r]]}$ by the map $\phi_T$. We claim that $C_{\tau}$ and $C_{\tau,T}$ are asymptotic. Indeed for all n we have:

$$ d(C_{\tau}[n],C_{\tau,T}[n])=d(\tau[n], T \cup \tau[n+r])$$
$$\leq \rho(\tau[n], T \cup \tau[n+r])=|\tau[n] \Delta (T \cup \tau[n+r])|=r+|T|-r=|T|$$.

\begin{lema}
The Tits angle between any two distinct special points is $\pi/2$.
\end{lema}

\noindent \emph{Proof.} Let $\tau, \sigma \in \Game$ be two distinct long snakes. We may assume $\tau < \sigma$ (in the obvious lexicographic order). By standard $CAT(0)$ considerations to prove that $\angle(\xi(\tau),\xi(\sigma))=\pi/2$ it is enough to find a positive quadrant and two representatives of these points which are axes.

Consider $T=\tau \cap \sigma \in Snakes$ and notice that $T=\tau(r)=\sigma(r)$ for some positive $r$. Consider the special orthant $F(T,\{i,i+1\},\tau[[r]], \sigma[[r]])$, where $(i,i+1)$ is the position of the unique free caret of $T$. The axes of this quadrant are the half-line geodesics $C_{\tau, T}$ and $C_{\sigma, T}$, so we are done. $\clubsuit$

\vspace{2mm}

The following fact is an important technicality.

\begin{lema}
For any distinct $\tau_1, \tau_2, \ldots, \tau_n \in \Game$ there is $T \in Trees$ and a special $n$-dimensional orthant with origin in $T$ and axes $C_{\tau_1,T},C_{\tau_2,T}, \ldots, C_{\tau_n,T}$.
\end{lema}

\noindent \emph{Proof.} We may assume $\tau_1 < \tau_2 < \ldots < \tau_n$. There is a large positive $r$ such that for any $i \in [n]$, $\tau_i[r] \notin \tau_j$ for any $j \neq i$. Consider $T= \bigcup_{i=1}^n \tau_i[r] \in Trees$ and $I=\{i_1,i_2, \ldots, i_n\} \subseteq [|T|+1]$ with $i_1 < i_2 < \ldots < i_n$, where $i_k$  is the index of the leaf of $T$ where $\tau_k[r]$ is connected with $\tau_k[r+1]$ in $\tau_k$. Consider the $n$-dimensional orthant $F(T,I,\tau_1[[r]], \tau_2[[r]], \ldots, \tau_n[[r]])$. The axes of this orthant are $C_{\tau_1,T},C_{\tau_2,T}, \ldots, C_{\tau_n,T}$. $\clubsuit$

\vspace{2mm}

Along the lines of the previous argument, we make a remark about the axes of an arbitrary special orthant. Let $F(T,I,\sigma_1, \ldots, \sigma_n)$ be an $n$-dimensional orthant, where $I=\{i_1, \ldots, i_n\}$ with $i_1 < \ldots <i_n$. For any $k \in \{1, \ldots, n\}$, let $\tau_k$ be the unique long snake contained in $T \oplus^{i_k} \sigma_k$. The axes of $F(T,I,\sigma_1, \ldots, \sigma_n)$ are $C_{\tau_1,T}, \ldots, C_{\tau_n,T}$, so to any special orthant we can attach a finite set of special points on the boundary (or in other words a finite profile). Notice that $\tau_1 \leq \ldots \leq \tau_n$ and from now on when we talk about axes of a special orthant we will always enumerate them in increasing order.

\vspace{2mm}

Denote by $\mathcal{S}$ the set of all the points on the $CAT(0)$ boundary of $\mathcal{T}$ which correspond to all the half-line geodesics in all the special orthants. In order to index conveniently these points, we proceed in several steps. The half-line geodesics starting from the origin in a special orthant are in bijection with the points of positive coordinates on the corresponding unit sphere (by just taking the intersection with the sphere).

\vspace{2mm}

\emph{Step 1.} In a first instance, for any  special orthant $F$ (we ignore the other terms of the notations since they won't be relevant) with axis $\tau_1, \ldots, \tau_n$, we'll denote by $$\xi(F)(\tau_1, \ldots, \tau_n)(t_1, \ldots, t_n)$$
 the point on the $CAT(0)$ boundary corresponding to the half-line geodesic in $F$ starting from origin and crossing the unit sphere in $F$ at coordinates $(t_1, \ldots, t_n)$. Notice that $t_i \geq 0$ for all $i=1, \ldots, n$ and $\sum_{i=1}^nt_i^2=1$.

\vspace{2mm}

\emph{ Step 2.} If we have two orthants of the same dimension with  pairwise asymptotic axes (in order) then the half-line geodesics at the same coordinates in the two flats are also asymptotic. In particular, if we have two special orthants $F_1$ and $F_2$ with the same axes $\tau_1, \ldots, \tau_n$, then we have $$\xi(F_1)(\tau_1, \ldots, \tau_n)(t_1, \ldots , t_n)=\xi(F_2)(\tau_1, \ldots, \tau_n)(t_1, \ldots, t_n)$$ So actually at this point we can drop some notations and simply write $\xi(\tau_1, \ldots, \tau_n)(t_1, \ldots, t_n)$ for the point corresponding to geodesics with prescribed coordinates on prescribed axis, independent of the orthant chosen with these properties.

\vspace{2mm}

\emph{Step 3.} In the last notation from the previous step it is enough to consider only the positive coordinates and the corresponding points, since we have $\xi(\tau_1, \ldots , \tau_n)(t_1, \ldots, t_n)= \xi(\tau_{i_1}, \ldots, \tau_{i_k})(t_{i_1}, \ldots, t_{i_k})$, where $\{i_1, \ldots, i_k\} \subseteq [n]$ is the set of all the indices $j$ with the property that $t_j >0$. Indeed, if we consider the geodesic representing the point $\xi(\tau_1, \ldots , \tau_n)(t_1, \ldots, t_n)$ in some appropriate orthant, we can actually consider the special "sub-orthant" generated only by the axes corresponding to non-zero coordinates. So we can drop the zero coordinates.

\vspace{2mm}

\emph{Step 4.} We have $\xi(\tau_1, \ldots, \tau_n)(t_1, \ldots, t_n)=\xi(\sigma_1, \ldots, \sigma_m)(s_1, \ldots, s_m)$ if and only if $m=n, \tau_i=\sigma_i, t_i=s_i$ for $i=1, \ldots, n$. We can construct special orthants with prescribed points corresponding to the axes: we consider a special positive flat with axis $\{\tau_1, \ldots, \tau_n\} \cup \{\sigma_1, \ldots, \sigma_m\}$. We can choose the geodesics corresponding to our two points and the conclusion becomes obvious.

\vspace{2mm}

After these reductions, we can write $\mathcal{S}=\{ \xi(\tau_1, \ldots, \tau_n)(t_1, \ldots, t_n) | \tau_1 < \ldots < \tau_n, \tau_i \in \Game, n \geq 1, t_i >0, \sum_{i=1}^n t_i^2=1 \}$. Let $l_{0,1,+}^2(\Game)=\{\phi \in l^2(\Game) | \|\phi\|=1,\phi \geq 0,supp(\phi)< \infty \}$. There is a simple bijection between this set and $\mathcal{S}$: to each $(\tau_1, \ldots, \tau_n)(t_1, \ldots, t_n)$ we can attach an element $\phi \in l_{0,1,+}^2(\Game)$, defined by $\phi(\tau_i)=t_i$ for $i=1, \ldots, n$ and $\phi(\sigma)=0$ if $\sigma \notin \{\tau_1, \ldots, \tau_n \}$. So from now on we will write $\xi(\phi)$ with $\phi \in l_{0,1,+}^2(\Game)$ for the elements of $\mathcal{S}$.

\begin{lema}
The Tits angle between $\xi(\phi)$ and $\xi(\psi)$ is $arccos(<\phi,\psi>)$ for $\phi, \psi \in l_{0,1,+}^2(\Game)$ and $<,>$ is the scalar product in $l^2(\Game)$.
\end{lema}

\vspace{2mm}

\emph{Proof.} To compute the angle between $\xi(\phi)$ and $\xi(\psi)$ it is enough to look at special orthant  that contains representatives of both points. Let $F$ be a special orthant with axis representing $supp(\phi) \cup supp(\psi)$. Let $(t_1, \ldots, t_n)$ be the coordinates corresponding to $\xi(\phi)$ and $(s_1, \ldots, s_n)$ be the coordinates corresponding to $\xi(\psi)$. We have that $\angle(\xi(\psi),\xi(\phi)=arccos(\sum_{i=1}^n t_is_i)=arccos(<\phi,\psi>)$ (the last inequality follows from the fact that zero coordinates that we added do not influence the scalar product). $\clubsuit$

\vspace{2mm}

By density we have Corollary 2, that the Tits boundary of $\mathcal{\tau}$ contains the positive part of the unit sphere of $l^2(\Game)$. In section 3.3. we put this construction in a much larger, but cleaner perspective.

\subsection{The Christmas Tree}

In this section we give some intuition for how the rest of the boundary looks like. Developments follow shortly  in the next section. We construct a point on the boundary which is not covered by the previous construction and it will be quite clear that many similar points exist.

Recall that $T_k$ is the full (finite) tree of depth $k$ and  consider the curve starting at origin and going (diagonally) at the geodesic speed to $T_1$ and than in the same way to $T_2$, $T_3$... One can think about this curve as going as diagonally possible through "the center" of the complex. We call this curve \emph{the Christmas Tree} since intuitively covers the full binary trees which look like a Christmas trees. Precisely, $\chi:[0,\infty) \rightarrow Trees$, determined by $\chi(0)=*, \chi(1)=T_1, \chi(1+\sqrt{2})=T_2,..., \chi(1+\sqrt{2}+\sqrt{4}+...+ \sqrt{2^{n-1}})=T_n$ for any $n$ (and between these points we take the corresponding geodesics).

\begin{lema}
$\chi$ is a geodesic half-line.
\end{lema}

\noindent \emph{Proof.} We only have to check that $\chi$ is a local geodesic near the points $T_k$'s. Fix a positive integer $k \geq 1$ and let $t > 0$ such that $\chi(t)=T_k$. Let $\varepsilon >0$ be small enough. The point $c(t-\varepsilon)$ will has the reduced diagram $T_{k-1}[1-\frac{\varepsilon}{2^{(k-1)/2}},\ldots, 1-\frac{\varepsilon}{2^{(k-1)/2}}]$ and $c(t+\varepsilon)$ has the reduced form $T_k[\frac{\varepsilon}{2^{k/2}}, \ldots, \frac{\varepsilon}{2^{k/2}}]$. By \ref{calcul recursiv}, via the lantern $\phi_{T_{k-1}}$ we have

\begin{eqnarray*}
d(c(t+\varepsilon),c(t-\varepsilon))&=&
\sqrt{\sum_{i=1}^{2^{k-1}}d(*[1-\frac{\varepsilon}{2^{(k-1)/2}}],\wedge[\frac{\varepsilon}{2^{k/2}},\frac{\varepsilon}{2^{k/2}}])^2}\\
&=& \sqrt{2^{k-1}(\frac{\varepsilon}{2^{(k-1)/2}}+\sqrt{\frac{\varepsilon^2}{2^k}+\frac{\varepsilon^2}{2^k} })^2}=2\varepsilon.
\end{eqnarray*}
$\clubsuit$

It is clear that the profile of this geodesic is the full binary tree, as it is also clear that the profiles of the points corresponding to the special orthants from the previous section correspond to the finite profiles, that is closed trees which are finite union of finitely many long snakes. It is also quite intuitive that a Christmas-like construction can be performed for the wide range of intermediate profiles between the finite ones and the full one. As we said, we will put the Christmas Tree in context, but for the moment let us just record its quite mysterious presence.

\subsection{Drawing geodesics}

Warning: the next paragraph describes informally  the process of building the half-line geodesics. Since the text might be confusing (for we avoid some technical issues), the reader can skip directly to the rigorous mathematics which is done after.

In this section we describe all the half-line geodesics of $\mathcal{T}$ which starts at the origin. Our method is very simple. A half-line geodesic is determined by the configuration of its linear pieces on each maximal cube it crosses. In turn, each piece is determined by its angles with the axes of the cube to which belongs: each geodesic $c$ starts at the origin and we have $c(1)=\wedge$, then we reach the cube $C(\wedge)$; here once we fix the angles of $c$ with the axis $[\wedge, L_2]$ and $[\wedge,R_2]$ the geodesics it is determined until it hits a face with a new (maximal) cube ; there we have to check the angles with the axes again... the correct "passing" conditions on angles will be: at the axes which are common for the old cube the angle is conserved, but at the new axes (which will be determined by a new caret occurring in the reduced form of the geodesic at that time) the passing condition will be $cos^2(leftson)+ cos^2(rightson)=cos^2(theirparent)$, where $cos(leftson)$ means the $cosinus$ of the angle of the geodesic with the axes determined by the left son etc. As one advances with describing the geodesic, by looking at the reduced diagram of the geodesic (from time to time), will notice that the diagrams "draw" a larger and larger tree, which at "the end" of the process will be a closed tree, representing the profile of the geodesic.

\vspace{2mm}

We now move on to facts. We need some terminology.

\begin{definitie}
\label{desene}

(a) If $x \in \mathcal{T}$, if $T[t_1, \ldots, t_s]$ is its reduced diagram, then we call $T$ its $\emph{mother}$ and the maximal cube $C(T)$ is called \emph{its mother cube}. We denote by $[x]$ the mother of $x$.

(b) An infinite sequence of distinct trees $S_1 \subseteq S_2 \subseteq \ldots$ is called \emph{a discrete drawing} if

(i) $S_1=*$

(ii) $S_k \subset S_{k+1} \subset S_k \oplus$

(iii) $\bigcup_{k \geq 1}S_k$ is a closed tree, called \emph{the final drawing}.

(c) \emph{A continuous drawing} or \emph{a flow} will be a map $V$ from the vertices of the full binary tree to $[0,1]$ with the property that $V$ maps the root to $1$ and the value of $V$ on a vertex is the sum of the values on its children. If we identify the vertices of the full binary tree with standard dyadic intervals (as usual) the condition reads:

(i) $V([0,1])=1$

(ii) $V([\frac{k}{2^n},\frac{k+1}{2^n}])=V([\frac{2k}{2^{n+1}},\frac{2k+1}{2^{n+1}}])+ V([\frac{2k+1}{2^{n+1}},\frac{2k+2}{2^{n+1}}])$

(d) If $c$ is a geodesic starting at the origin, we call $c(t)$ for some $t \geq 0$ a \emph{corner point} of $c$ if there is an $\varepsilon$ such that $[c(t-\varepsilon),c(t)]$ belongs to  a same cube and $[c(t),c(t+\varepsilon(t)]$  belongs a same cube, but there is no cube to which  $c[t-\varepsilon, t+\varepsilon]$ belongs. The number (moment of time) $t$ is called a \emph{corner time/moment}. We slightly violate the definition to allow $c(0)=*$ to be a corner point and $0$ a corner moment.
\end{definitie}

Warning: next, we regard a tree as embedded in the full binary tree.

To each flow $V$ we will associate a half-line geodesic $c_V$ starting at the origin in $\mathcal{T}$. We construct by induction a sequence of points $(x_n)_{n \geq 0}$ in $\mathcal{T}$ which will turn out to be the corner points of $c_V$. We set $x_0=*$ and $x_1=\wedge$. Before giving the general formula we construct one more point. We define $x_2$ to be the unique point in the mother cube of $x_1$ such that the mother of $x_2$ is different from the mother of $x_1=\wedge$ and the angles of $[x_1,x_2]$ with the axes of the cube $C(\wedge)$ are taken from the flow $V$ as follows:

$$\angle_{x_1}([x_1,x_2],[\wedge,L_2])=arccos(\sqrt{V([0,1/2])}$$

and

$$\angle_{x_1}([x_1,x_2],[\wedge,R_2])=arccos(\sqrt{V([1/2,1])}.$$

\noindent Notice that by construction $\wedge \subseteq [x_2] \subseteq T_2$.

We continue by induction and define $x_{n+1}$ to be the unique point in the mother cube of $x_n$ with the following properties:

\vspace{2mm}

\noindent - $[x_n] \neq [x_{n+1}]$

\vspace{2mm}

\noindent - the angles made by the geodesic segment $[x_n,x_{n+1}]$ with the axes of $C([x_n])$, namely  with $[[x_n],[x_n]\oplus^j]$ for $j=1, \ldots, |[x_n]|+1$, are given by $\alpha_j=arccos(\sqrt{V(P_j)})$, where $P_j$ is the $j^{th}$ leaf of $[x_n]$ (viewed as a vertex in the full binary tree). Of course, $\sum_{j}cos^2(\alpha_j)=1$.

\vspace{2mm}

Notice that, by construction, we have $[x_n] \subseteq [x_{n+1}] \subseteq [x_n]\oplus$ (with the first inclusion always strict).

We consider the following sequence of non-negative real numbers (will be the corner moments of $c_V$): $t_0=0$, $t_1=1$ and $t_{n+1}=t_n+d(x_n,x_{n+1})$. By construction, since the sequence of mothers is strictly increasing, we have $t_n \rightarrow \infty$. We define $c_V(t_n)=x_n$ for any $n \geq 0$ and extend $c_V$ on each interval $[t_n,t_{n+1}]$ with the geodesic between $x_n$ and $x_{n+1}$ (which travels in the cube $C([x_n])$- the entering point is $x_n$ and the exit point is $x_{n+1}$). We notice that the map $c_V:[0, \infty) \rightarrow \mathcal{T}$ is increasing, that is, $c_V(t) \subseteq c_V(s)$ if $t \leq s$. We also notice that by construction, the sequence of mothers $[c_V(t_0)] \subseteq [c_V(t_1)] \subseteq \ldots$ forms a discrete drawing with the final drawing the closed tree obtained by from the full tree by deleting all the subtrees rooted at those $P$ with $V(P)=0$.

We first do some examples and then show that $c_V$ is a geodesic ray.

\begin{figure}[!h]
\label{rolex}
\begin{center}
\includegraphics{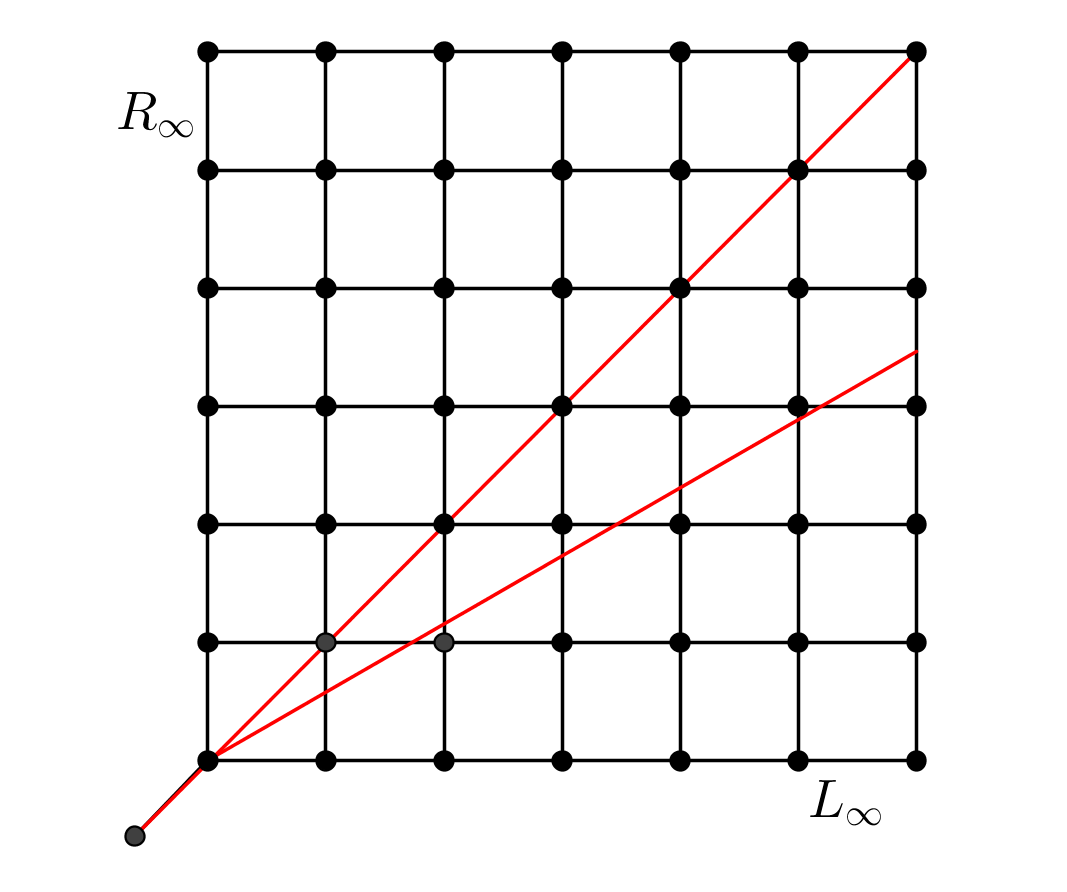}
\caption{In the figure we have two geodesic rays with the images included in "the quadrant with tail" generated by $L_{\infty}$ and $R_{\infty}$. The diagonal one correspond to the flow defined by $V([0,2^{-n}])=1/2$ and $V([1-2^{-n},1])=1/2$ for any $n \geq 1$. The other one is "more generic", it correspond to the flow $V([0,2^{-n}])=\alpha$ and $V([1-2^{-n},1])=1-\alpha$ for any $n \geq 1$, if the angle of the geodesic with $L_{\infty}$ is $arccos{\sqrt{\alpha}}$. Both have as final drawing the tree $L_{\infty} \cup R_{\infty}$, but the discrete drawing of the diagonal one is $* \subseteq \wedge \subseteq T_2 \subseteq T_3 \ldots$, while for the other one is $* \subseteq \wedge \subseteq L_2 \subseteq T_2 \subseteq L_3 \cup R_2 \subseteq L_4 \cup R_2 \subseteq L_4 \cup R_3 \subseteq L_5 \cup R_3 \ldots$.}
\end{center}
\end{figure}

\begin{exemplu}

-If $V$ is the flow which maps each standard dyadic interval to its Lebesque measure, then $x_n=T_n$ for any $n \geq 0$ (recall that $T_n$ is the full tree of depth $n$) and $c_V$ is the Christmas Tree. The final drawing is the full binary tree.

-If $V$ is the flow which maps each interval of the form $[0,2^{-n}]$ to $1$ and all the other standard dyadic intervals to zero, then $x_n=L_n$ for any $n$ and $c_V=c_{L_{\infty}}$. The final drawing is $L_{\infty}$.

-More generally, for any snake $\sigma$ the associated geodesic $c_{\sigma}$ is defined by the flow mapping all the intervals corresponding to the interior vertices of $\sigma$ to $1$ and all the others to zero. The final drawing is $\sigma$.

-If $V$ is the flow which maps all the intervals of the form $[0,2^{-n}]$ and $[1-1/2^n,1]$ to $1/2$ and all the others intervals to zero, then $x_n=L_{n} \cup R_{n}$ and $c_V$ is the diagonal of the two-dimensional orthant generated by $L_{\infty}$ and $R_{\infty}$.

\end{exemplu}

\begin{lema}
\label{geodezice}
For any flow $V$, the curve $c_V$ is a half-line geodesic starting at the origin. Its profile is the closed tree obtained from the full binary tree by deleting all the subtrees rooted at vertices $P$ with $V(P)=0$.
\end{lema}

\noindent \emph{Proof.} Let $V$ be a flow and let $c=c_V$ be the corresponding curve as defined above. We only have to check that $c$ is a local geodesic near the points $x_n$ (by construction the $x_n$'s will be the corner points of the geodesic). Let $t=t_n >0$ be an element in the sequence of moments in the construction of $c_V$, in particular $c(t)=x_n$. Fix a very small $\varepsilon >0$. Denote by $S$ the mother of $c(t-\varepsilon)$ and by $T=S \oplus^I$ the mother of $c(t)$ (and in the same time the mother of $c(t+ \varepsilon)$). Denote by $\alpha_j$ the angles between $c$ and $[S,S\oplus^j]$ in the cube $C(S)$. If $j \notin I$ then the angle in the cube $C(T)$ between $c$ and $[T, T \oplus^{\sigma(j)}]$ is also $\alpha_j$ by construction, where $\sigma(j)$ denotes the label of the $j^{th}$ leaf of $S$ when viewed in $T$. If $i \in I$ and if we denote by $\alpha_{i}^l$ and $\alpha_{i}^r$ the angles in $C(T)$ between $c$ and the axes determined by the left and right leaf of the free caret of $T$ which is attached at the $i^{th}$ leaf of $S$ (when viewed in $T$), then we have by construction $cos^2(\alpha_i^l)+cos^2(\alpha_i^r)=cos^2(\alpha_i)$.

Denote by $r_j$ the coordinate of $c(t-\varepsilon)$ at the $j^{th}$ leaf. If $j \notin I$ we denote by $t_j$ the coordinate of $c(t)$ at the $j^{th}$ leaf of $S$ when viewed as a leaf in $T$ and by $s_j$ the coordinate of $c(t+\varepsilon)$ at the same position. If $i \in I$ then we denote by $t_i^l$ and $t_i^r$ the coordinates of $c(t)$ at the left and right children of the $i^{th}$ leaf of $S$ when viewed as a vertex in $T$; $s_i^l$ and $s_i^r$ denotes the coordinates of $c(t+ \varepsilon)$ at the same positions. By the construction we have $t_j-r_j=s_j-t_j=\varepsilon cos(\alpha_j)$ for $j \notin I$, and $1-r_i=\varepsilon cos(\alpha_i)$, $s_i^l=\varepsilon cos(\alpha_i^l)$, $s_i^r=\varepsilon cos(\alpha_i^r)$ for $i \in I$. Using \ref{calcul recursiv} with the lantern $\phi_S$ we have:

\begin{eqnarray*}
d(c(t+\varepsilon),c(t-\varepsilon))^2&=& \sum_{j \notin I} (s_j-r_j)^2 + \sum_{i \in I}[1-r_i+ \sqrt{(s_i^l)^2+(s_i^r)^2}]^2\\
       &=& \sum_{j \notin I} [(s_j-t_j)+(t_j-r_j)]^2 + \sum_{i \in I}[ (1-r_i)^2 + (s_i^l)^2 + (s_i^r)^2]\\
       &+&2\sum_{i \in I}(1-r_i)\sqrt{(s_i^l)^2+(s_i^r)^2}\\
       &=&\sum_{j \notin I}(t_j-r_j)^2 + \sum_{i \in I}(1-r_i)^2\\
       &+& \sum_{j \notin I}(s_j-t_j)^2 + \sum_{i \in I}[(s_i^l)^2+(s_i^r)^2]\\
       &+& 2\sum_{j \notin I}(s_j-t_j)(t_j-r_j) + 2\sum_{i \in I}(1-r_i)\sqrt{(s_i^l)^2 + (s_i^r)^2}\\
       &=& d(c(t),c(t-\varepsilon))^2 + d(c(t),c(t+\varepsilon))^2 \\
       &+& 2\varepsilon^2\sum_{j \notin I}cos^2\alpha_j + 2\varepsilon^2\sum_{i \in I}cos\alpha_i \sqrt{(cos\alpha_i^l)^2+(cos\alpha_i^r)^2}\\
       &=& 2\varepsilon^2+ 2\varepsilon^2(\sum_{j \notin I}cos^2\alpha_j + \sum_{i \in I}cos^2\alpha_i)\\
       &=& 4\varepsilon^2
\end{eqnarray*}

$\clubsuit$

\vspace{2mm}

We will show later that all the half-line geodesics starting at the origin in $\mathcal{T}$ are of the form $c_V$ for some flow $V$. We are now interested in describing the half-lines geodesics asymptotic with a given $c_V$, but starting at an arbitrary point $x \in \mathcal{T}$.

So we fix a flow $V$ and a point $x \in \mathcal{T}$. We define the geodesic $C_{V,x}$ in the same way we defined $c_V$; namely for each linear piece, starting with the one in the mother cube of $x$, we fix the angles with the axes of the mother cube by reading the number written at the vertex (leaf) which determines the axis and taking the $arccos$ of its square root. In this way we obtain again a monotone geodesic.

\begin{lema}
\label{asi}
For any flow $V$ and any $x \in \mathcal{T}$, $c_V$ and $c_{V,x}$ are asymptotic.
\end{lema}

Before doing the proof we present some examples:

\begin{exemplu}
(a) Consider $V$ the flow corresponding to $L_{\infty}$ and the point $x=L_5 \cup R_3$. If $c_1=c_V$ and $c_2=c_{V,x}$ then for any $t \geq 5$, if $s$ denotes the largest integer smaller than $t$, by \ref{calcul recursiv}, we have
$$d(c_1(t),c_2(t))=d((L_n[t-n,0, \ldots,0],(L_n \cup (R_2 \oplus^2 L_2))[n-t,0, \ldots,0]))=\sqrt{29}$$\\
(b) Consider $V$ the flow corresponding to The Christmas Tree and $x=L_2$. If $c_1=c_V$, $c_2=c_{V,x}$ and  $s_k=1+ \sqrt{2}+ \ldots + \sqrt{2^k}$, then for any large $n$, $c_2(s_n)$ will have mother $T_n$ and coordinates $t_k=(\sqrt{2}+1)2^{-n/2}$ for $k \leq 2^{n-1}$ and $t_k=2^{-n/2}$ for $k > 2^{n-1}$. A simple computation shows
$$d(c_1(s_n),c_2(s_n))=\sqrt{2+\sqrt{2}}$$
\end{exemplu}

\noindent \emph{Proof of the lemma.} An important fact to keep in mind: the intersection of the (images) geodesics $c_V$ and $c_{V,x}$ with the subtree (of the full binary tree) rooted at a vertex $P$ grows at geodesical speed $\sqrt{V(P)}$. Denote $c_1=C_V$ and $c_2=C_{V,x}$. Let $T$ be the mother of $x=T[s_1, \ldots, s_m]$. Denote by $J$ the set of labels of leaves of $T$ at which $V$ is non-zero. We call these leaves alive and all the others dead. For any alive leaf $j$, let $t_j$ be the moment of time when $c$ reaches the leaf $j$ at the coordinate $s_j$, that is, $Pr_{\sigma}(c(t_j))=Pr_{\sigma}(T)$ (\ref{proiectie serpi}), where $\sigma$ is a long snake containing the snake determined by the leaf $j$. Let $t_0=max\{t_j | j \in J\}$. Fix any $t \geq t_0$. Let $S$ be the mother of $c_1(t)$. For any leaf (labeled) $k$ of $S$ at which $V$ is not zero, we denote by $\sigma(k)$ the (label of the) unique  alive leaf of $T$ such that the $k^{th}$ leaf of $S$ is included in a subtree (of the full tree) rooted at the $\sigma(k)^{th}$ leaf of $S$ (the subtrees rooted at the dead leaves of $T$ are disjoint from the image of $c$, by construction).

If $P_k$ denotes the vertex at the $k^{th}$ leaf of $T$, then the geodesic $c_2$ reaches the coordinate of $c_1(t)$ at $P_k$ in $t-t_{\sigma(k)}$ time (and then it still has $t_{\sigma(k)}$ time to grow at geodesic speed $\sqrt{V(P_k)}$). By \ref{calcul recursiv} with the lantern $\phi_S$ and using that $\sum_kV(P_k)\leq 1$, we have:

$$d(c_1(t),c_2(t))^2= \sum_{k}t^2_{\sigma(k)}V(P_k) + const \leq t_0^2 + const$$

\noindent where $const$ is the constant quantity resulting from computing at the dead leaves of $T$ ( recall that the subtrees rooted at the dead leaves do not grow along $c_1$ and $c_2$); for example, using the median distance, $const \leq |T-\tau|$ where $\tau$ is the closed tree giving the profile of $c_1$. $\clubsuit$

\vspace{2mm}

For any flow $V$ denote by $\xi_V$ the point at infinity represented by $c_V$. $I_{n,k}$ denotes the standard dyadic interval $[\frac{k-1}{2^n},\frac{k}{2^n}]$.

\begin{lema}
If $V$ and $W$ are two flows then the Tits angle between $\xi_V$ and $\xi_W$ is given by the formula:
$$arccos( lim_{n \rightarrow \infty}\sum_{k=1}^{2^n}\sqrt{V(I_{n,k})W(I_{n,k})})$$
\end{lema}

\noindent \emph{Proof.} Let $x$ be any point in $\mathcal{T}$. If $S$ is the mother of $x$, then the cosine of the angle at $x$ between $c_{V,x}$ and $c_{W,x}$ equals $Y_S=\sum_j \sqrt{V(P_j)W(P_j)}$, where $j$ runs over all the leaves of $S$ and $P_j$ is the vertex in the full tree at the $j^{th}$ leaf of $S$. These angles depend only on mothers, so in order to compute the cosine of the Tits angle we only have to find $inf\{Y_S | S \in Trees\}$. By the mean inequality we have $Y_S \leq Y_T$ if $ T \subseteq S$. In particular, the infimum we search for can be computed as the limit of the decreasing sequence $(y_{T_n})$ (or along any increasing sequence of trees with the union the full binary tree). $\clubsuit$

\vspace{2mm}

To finish the proof of Theorem3, we still  have to deduce the formula for the action of $F$ and to show that the points associated with the flows are the only ones in the boundary of $\mathcal{T}$.

\begin{lema}
Any half-line geodesic in $\mathcal{T}$ starting at the origin is of the form $c_V$ for some flow $V$.
\end{lema}

\noindent \emph{Proof.} Let $c$ be a geodesic half-line starting at the origin. Let $t_1=0,t_2, \ldots,$ be the corner moments.

\vspace{2mm}

\emph{Step1.} $c$ is monotone, that is, for any $s \leq t$ we have $c(s) \leq c(t)$; in particular $[c(t_1)] \subseteq [c(t_2)] \subseteq \ldots$ form a discrete drawing.

\vspace{2mm}

We show this by induction, on each linear piece $[c(t_i),c(t_{i+1})]$. The first step is trivial. Suppose now that we have shown that $c$ is monotone on the interval $[0,t_n]$. Fix a very small $\varepsilon > 0$ and also denote for simplicity $t=t_n$. Next we analyze the points $c(t-\varepsilon), c(t), c(t+ \varepsilon)$. The notations are similar with the ones in the proof of Lemma \ref{geodezice}. Let $S$ be the mother of $c(t-\varepsilon)$ and $T$ the mother of $c(t)$ such that $T=S \oplus^I$ for some subset $I$ of the leaves of $S$. We denote by $J$ the set of all indices of the leaves of $S$ which are not in $I$. The difference now is that the mother of $c(t+\varepsilon)$ can be apriori different from $T$ (if and only if the geodesic is not monotone). We need to further partition the sets $I$ and $J$. Let $I_1$ be the subset of $I$ consisting of all the $i \in I$ such that the caret at the $i^{th}$ leaf of $S$ doesn't belong to the mother of $c(t+\varepsilon)$. $I_2$ denotes all the rest of $I$, hence $I=I_1 \sqcup I_2$. We denote, like in the proof of \ref{geodezice}, $r_i$ and $s_i$ the coordinates of $c(t-\varepsilon)$ and $c(t+\varepsilon)$  for $i \in I_1$ and $r_i$ and $s_i^l,s_i^r$ the coordinates of $c(t-\varepsilon)$ and $c(t+ \varepsilon)$ for $i \in I_2$. Also, denote $s_j,t_j,r_j$ the coordinates of $c(t-\varepsilon),c(t),c(t+\varepsilon)$ for $j \in J$. We consider the following partition of $J$. $J_1$ denotes the set of all those $j \in J$ for which $t_j < s_j$ and $J_2$ all the others elements of $J$. Hence $J=J_1 \sqcup J_2$. Since $\varepsilon >0$ is very small, these partitions cover all the possibilities for the behavior of $c$ right after the moment $t$. $c$ is monotone after $t$ if and only if $I_1=J_1=\emptyset$. We show next that if one of the sets $I_1, J_1$ is non-empty, then $c$ is not locally geodesic near $t$. So we assume this.

Using \ref{calcul recursiv} with the lantern $\phi_S$, we have:

$$\varepsilon^2=d(c(t-\varepsilon,c(t))^2=\sum_{i \in I}(1-r_i)^2 + \sum_{j \in J}(t_j-r_j)^2$$

and

\begin{eqnarray*}
\varepsilon^2=d(c(t),c(t+ \varepsilon))^2&=& \sum_{i \in I_1}(1-s_i)^2 + \sum_{i \in I_2}[(s_i^l)^2+(s_i^r)^2]+\\
&+& \sum_{j \in J_1} (s_j-t_j)^2 + \sum_{j \in J_2}(t_j-s_j)^2
\end{eqnarray*}

We can now give the main computation:

\begin{eqnarray*}
d(c(t-\varepsilon),c(t+\varepsilon))^2 &=& \sum_{i \in I_1}(s_i-r_i)^2 + \sum_{i \in I_2}(1-r_i+ \sqrt{(s_i^l)^2+(s_i^r)^2})^2\\
&+& \sum_{j \in J_1} (s_j-r_j)^2 + \sum_{j \in J_2}(s_j-r_j)^2\\
&=& \sum_{i \in I_1}[(s_i-1) + (1-r_i)]^2 + \sum_{i \in I_2}(1-r_i+ \sqrt{(s_i^l)^2+(s_i^r)^2})^2\\
&+& \sum_{j \in J_1}[(t_j-r_j)+(s_j-t_j)]^2 + \sum_{j \in J_2}[(s_j-t_j)+(t_j-r_j)]^2\\
&=& d(c(t-\varepsilon),c(t))^2 + d(c(t),c(t+ \varepsilon))^2 \\
&+& 2\sum_{i \in I_2}(1-r_i)\sqrt{(s_i^l)^2+(s_i^r)^2} + 2\sum_{j \in J_2}(t_j-r_j)(s_j-t_j)\\
&-& 2\sum_{i \in I_1}(1-s_i)(1-r_i)-2\sum_{j \in J_1}(t_j-r_j)(t_j-s_j)\\
&\leq& 2\varepsilon^2 + \sum_{i \in I_2}(1-r_i)^2+\sum_{i \in I_2}[(s_i^l)^2+(s_i^r)^2] \\
&+& \sum_{j \in J_2}(t_j-r_j)^2+ \sum_{j \in J_2}(s_j-t_j)^2\\
&-& 2\sum_{i \in I_1}(1-s_i)(1-r_i)-2\sum_{j \in J_1}(t_j-r_j)(t_j-s_j)\\
&=& 4\varepsilon^2 - 2\sum_{i \in I_1}(1-s_i)(1-r_i)-2\sum_{j \in J_1}(t_j-r_j)(t_j-s_j)\\
&<& 4\varepsilon^2
\end{eqnarray*}

The last inequality follows if $I_1 \neq \emptyset$ or $J_1 \neq \emptyset$ and just before we applied the mean inequality for each term of the sums having coefficient +2 in front. So $c$ must be monotone to stand a chance of being a geodesic.

\vspace{2mm}

\noindent \emph{Step 2.} $c=c_V$ for some flow $V$.

\vspace{2mm}

Let $t$ be any corner moment and let $\varepsilon > 0$ be very small. Let $S$ be the mother of $c(t-\varepsilon)$ and let $T=S \oplus^I$ be the mother of $c(t)$, for some subset $I$ of the set of leaves of $S$ (the fact that $T$ has such a form follows from Step 1). $T$ is also the mother of $c(t+\varepsilon)$. Similarly with the proof of \ref{geodezice}, for the $j^{th}$ leaf of $S$ let $\alpha_j$ be the angle of $c$ with the axes of $C(S)$ determined by the $j^{th}$ leaf. For any $i \in I$ we denote by $\beta_i^l$ and $\beta_i^r$ the angles between $c$ and the axes of $C(T)$ determined by the left and right son of the $i^{th}$ leaf of $S$ when viewed in $T$. For any $j \notin I$, denote by $\beta_j$ the angle between $c$ and the axis of $C(T)$ determined by the $j^{th}$ leaf of $S$ when viewed in $T$. To show that $c$ is obtained from a flow it is enough to show that $\alpha_j=\beta_j$ for any $j \notin I$ and $cos^2(\alpha_i)=cos^2(\beta_i^l)+cos^2(\beta_i^r)$ for any $i \in I$. The underlying flow will have at a given vertex the cosine squared of the  value of the unique angle associated with the given vertex.

We have to consider again the notations for $c(t-\varepsilon),c(t), c(t+\varepsilon)$ in Step1, with the difference that now we can assume $I_1=J_1=\emptyset$ and hence $I_2=I$, $J_2=J$. In the new, much simpler setting, the computation in Step1 becomes:

\begin{eqnarray*}
d(c(t-\varepsilon), c(t+\varepsilon))^2 &=& d(c(t-\varepsilon),c(t))^2 + d(c(t), c(t+ \varepsilon))^2 \\
                                         &+& 2\sum_{i\in I}(1-r_i)\sqrt{(s_i^l)^2+(s_i^r)^2} + 2\sum_{j \in J}(t_j-r_j)(s_j-t_j)
                                         \end{eqnarray*}

By applying the mean inequality to each term of the two sums above, we get that $d(c(t-\varepsilon), c(t+\varepsilon))^2 \leq 4\varepsilon^2$ with equality if and only if for any $i \in I$ we have $1-r_i=\sqrt{(s_i^l)^2 + (s_i^r)^2}$ which implies $cos^2(\alpha_i)=cos^2(\beta_i^l)+cos^2(\beta_i^r)$, and for any $j \notin I$ we must have $s_j-t_j=t_j-r_j$ which implies $cos(\alpha_j)=cos(\beta_j)$. $\clubsuit$

\vspace{2mm}

We now compute the action of $F$ on the boundary of $\mathcal{T}$. Notice that a flow is determined if we specify the values at all but finitely many standard intervals. Also, notice that if $g \in F$, then $g$ maps linearly a standard interval to a standard interval with finitely many exceptions (namely the intervals corresponding to the inner vertices of the left tree of its reduced diagram). Now we can specify the action. For any $g \in F$ and any flow $V$ we have $(gV)(P)=V(g(P))$, where $P$ is any standard interval which is mapped by $g$ linearly to another standard interval. To see this action in term of diagrams, one has to take large trees $T$ such that the graphical composition with $g=(S,T)$ can be performed. The result is that the interval corresponding to the $i^{th}$ leaf of $S$ should have the value of $V$ at the interval corresponding to the $i^{th}$ vertex of $T$.

\begin{lema}
The boundary of $\mathcal{T}$ is invariant under $F$ and the action is given by $(gV)(P)=V(g(P))$ for any $g \in F$, any flow $V$ and any standard interval $P$ which is mapped linearly by $g$ onto another standard interval.
\end{lema}

\noindent \emph{Proof.} Consider $g \in F$ and $V$ a flow. First, for simplicity, we consider $V$ has the full binary tree as its final drawing. Let $S_1=* \subseteq S_2 \subseteq \ldots$ be the discrete drawing of $c_V$ and denote $t_1, t_2, \ldots, $ the associated corner moments. There exists a positive integer $N$ such that for all $n \geq N$ $g$ has a diagram of the form $(R_n,S_n)$. Each cube $C(S_n)$ is mapped by $g$ onto the cube $C(R_n)$ and $g$ preserves the angles formed by the geodesic $c$ with the axes of $C(S_n)$. It follows that $gc_V(t)=c_{gV,gc(t_N)}(t-t_N)$ for any $t \geq t_N$. By \ref{asi} $gc$ represent the point at infinity we called $gV$ and we are done in this case.

Consider now the general case, when the final drawing of $V$ is a closed tree $\mathcal{\tau}$. We keep the notations of the previous case for the discrete drawing and corner moments of $c_V$. Denote by $(R,S)$ the reduced diagram of $g$. There is a positive integer $N$ such that $\tau \cap S \subseteq  S_N$ and hence for any $n \geq N$ there is a diagram $(\widehat{R}_n,S_n \cup S)$ representing $g$. $g$ sends the cubes $C(S_n)$ (for $n \geq N$) to the cubes $C(R_n, \Delta)$, where $\Delta$ (by slight abuse of notation) denotes only the configuration of non-trivial trees on the right, which we are forced to add in order to perform the multiplication. This configuration is the same for any $n \geq N$ and its total number of carets is bounded by $|S-\tau|$. Consider now for any $n \geq N$ the trees $R_n$ obtained from $\widehat{R}_n$ by deleting all the carets corresponding to $\Delta$, in other words keeping only the vertices of $\widehat{R_n}$, where we have non-zero values of $gV$. The geodesic $gc$ is asymptotic with the geodesic $c_{gV}$, which drawing is determined by the sequence of $R_n$'s. $\clubsuit$

\begin{exemplu}
For a better understanding of what is going on in the last paragraph of the the proof, we do an example. Let $V$ be the unique flow determined by $L_{\infty}$ and let $g=(T_2, L_2\oplus^2) \in F$. For any $n \geq 2$, we have $gL_n=(L \oplus^{-1}, \Delta)$, where $\Delta$ is the sequence of trees consisting only of empty trees $*$ with the exception of the second last position where we have a $\wedge$; the superscript $^{-1}$ denotes that we attach a caret at the last position. For any $n \leq 2$, we have $d(gL_n,L_n) \leq 2$.
\end{exemplu}

We now observe the fixed points of $F$ in this region of the boundary. Notice first that $F$ acts on the set of profiles and fixes exactly four of them, corresponding to the following closed trees: $L_{\infty}$, $R_{\infty}$, $L_{\infty} \cup R_{\infty}$ and to the full binary tree $Bin$ (comparing with Theorem 5.1 in \cite{farley08}, these are the profiles $\Delta_L,\Delta_R, \Delta_{L-R}$ and $\Delta_{\infty}$ in the notation there). It is true that we work in a subcomplex, but it follows from Section 4.1 in $\cite{farley08}$ (especially see Proposition 4.3 and Lemma 4.5 there) that in fact the geodesics with these profiles cannot leave the subcomplex $\mathcal{T}$; and indeed, as soon as a geodesic leaves this subcomplex will cross hyperplanes of different nature- we don't carry here the classification of all profiles in the full complex, since this was done in great detail by Farley in his paper. With these observations we in fact capture all the points at infinity having these profiles. It follows immediately from our formula that $F$ fixes all the geodesics in the profiles $L_{\infty}$, $R_{\infty}$ and $L_{\infty} \cup R_{\infty}$ and we get so far a Tits arc of length $\pi/2$. (Notice that the profile $L_{\infty} \cup R_{\infty}$ contains all the flows $V_{\alpha, \beta}$ (for any $\alpha$ and $\beta$ such that $\alpha+ \beta=1$) defined by $V(P)=\alpha$ and $V(Q)=\beta$ for any $P$ situated on the left wing of the full tree and for any $Q$ situated on the right wing.)

Farley showed that any other profile doesn't contain global fixed points, with the exception of the (very large) profile corresponding to the full binary tree. With our formula we can rule it out. Indeed, let $V$ be a flow having this profile, that is, taking no zero value. We have that $V([0,1/4]) < V([0,1/2])$. Pick any $g \in F$ such that $g$ maps $[0,1/2]$ linearly onto $[0,1/4]$. With our formula we have $gV([0,1/2])=V(g[0,1/2])=V([0,1/4]) < V([0,1/2])$, hence $V$ is not fixed by $F$. Corollary 3 is now proved.

\newpage

\appendix
\section{Empty Profiles}

We sketch a two-dimensional $CAT(0)$ cubical complex with exactly two profiles and one point in the boundary. This gives a counter-example to Conjecture 2.8 (1) in \cite{farley08}.

The idea is to take an approximation by cubulation of the following subset of the plane: $\{(x,y)| x,y \geq 0, y \leq \sqrt{x} \}$. More precisely, let $f:[0, \infty) \rightarrow [0,\infty)$ be the following discrete version of the sqrt function: $f(x)=[\sqrt{x}]$, where $[u]$ denotes the greatest integer smaller than $u$. Let $M=\{(x,y)| x,y \geq 0, y \leq f(x) \}$ with the shortest path $CAT(0)$ distance and with the obvious integer cubulation. We fix $(0,0)$ the origin.

There is an obvious geodesic half-line starting at the origin, the $Ox$ axis. Any other half-line geodesic would hit the graph of $f$ in a non-corner point after a while (due to sub-linearity of $f$) and from there it cannot be further extended to remain a geodesic (due to the concavity of $f$).

We split the hyperplanes of $M$ in two categories: horizontals and verticals. For any $n \geq 0$, let $H_n$ be the intersection of $M$ with the euclidian line $y=n+1/2$ and $V_n$ be the intersection of $M$ with the euclidian line $x=n+1/2$. These are all the hyperplanes of $M$.

The family of all $V_n^+$ is a profile and it is realized by the unique half-line $Ox$. Notice that something a bit counter-intuitive happens: the set of all $H_n$'s is not a profile since $V_n^+ \supseteq H_m^+$ for $m$ sufficiently larger than $n$ and hence the third axiom fails. Nevertheless the family of all the hyperplanes together form a profile, as easily can be seen.

\newpage
\begin{center}
\large \bf Conclusion
\end{center}

We conclude with two problems we consider more or less in range after this thesis, especially the first one.

\vspace{3mm}

\noindent \emph{Problem1.} Classify the isometry type of each element in $F$.

\vspace{5mm}

\noindent \emph{Problem2.} Find a description of the full Tits Boundary of the complex, in the style of Theorem3.

\vspace{5mm}

Most of all, we hope that there will be interesting applications of $CAT(0)$ geometry to the Thompson group.

\newpage

\bibliographystyle{amsalpha}
\bibliography{ref}
\end{document}